\documentclass[final,onefignum,onetabnum]{siamonline190516}

\pdfoutput=1

\usepackage[section]{placeins} 

\usepackage{algorithmic}
\usepackage[caption=false]{subfig}
\ifpdf
  \DeclareGraphicsExtensions{.eps,.pdf,.png,.jpg}
\else
  \DeclareGraphicsExtensions{.eps}
\fi

\usepackage{amssymb}
\usepackage{amsmath}
\usepackage{amsfonts}
\usepackage{upgreek}
\usepackage{siunitx}
\usepackage{listings}
\usepackage{multirow}
\usepackage{multicol}
\usepackage{bigdelim}
\usepackage{color}
\usepackage{stmaryrd}
\usepackage{soul}
\usepackage[symbol]{footmisc}
\usepackage{booktabs}
\usepackage[section]{placeins} 
\usepackage[titletoc]{appendix}
\usepackage{calc}
\usepackage[capitalize]{cleveref}
\usepackage[sort&compress,square,numbers]{natbib}

\usepackage{graphicx}
\usepackage{graphics}
\usepackage{float}
\usepackage{wrapfig}
\usepackage[percent]{overpic}
\usepackage{varwidth}
\usepackage{epstopdf}
\usepackage{adjustbox}
\usepackage{tikz}
\usepackage{pgfplots}
\usetikzlibrary{arrows,matrix,positioning,fit}
\usetikzlibrary{shapes,positioning}
\usetikzlibrary{backgrounds}
\pgfplotsset{compat=1.14}
\usepgfplotslibrary{colorbrewer}
\usepgfplotslibrary{patchplots}
\usepgfplotslibrary[colorbrewer]
\usetikzlibrary{pgfplots.colorbrewer}
\usetikzlibrary[pgfplots.colorbrewer]
\usepgfplotslibrary{fillbetween}
\usepgfplotslibrary{units}
\usetikzlibrary{spy}
\usepackage{pgfplotstable}
\graphicspath{ {Figs/} }

\newcommand\hb[1]{\hat{\mathbf{#1}}}

\newcommand\tb[1]{\tilde{\mathbf{#1}}}

\newcommand\px[2]{\frac{\partial #1}{\partial {#2}}}
\newcommand\pxi[3]{\frac{\partial^{#1}#2}{\partial {#3}^{#1}}}

\newcommand\dx[2]{\frac{\mathrm{d} #1}{\mathrm{d} #2}}

\newcommand{\half}{\frac{1}{2}}
\newcommand{\diag}{\text{diag}}

\newlength\myheight
\newlength\mydepth
\settototalheight\myheight{Xygp}
\newsiamremark{remark}{Remark}
\newsiamremark{hypothesis}{Hypothesis}
\crefname{hypothesis}{Hypothesis}{Hypotheses}
\newsiamthm{claim}{Claim}

\begin{document}

\title{On the Use of RBF Interpolation for Flux Reconstruction\thanks{\corresponding{Rob Watson(\email{r.a.watson@lboro.ac.uk})} \funding{None to declare}}}
\headers{RBF-FR}{R. Watson and W. Trojak}

\author{Rob Watson\thanks{Department of Aeronautical {\&} Automotive Engineering Aerospace Engineering, Loughborough University, Loughborough, LE11 3TU, UK} \and
Will Trojak\thanks{Department of Aeronautics, Imperial College London, South Kensington, London, SW7 2AZ}}

\maketitle

\begin{abstract}
Flux reconstruction provides a framework for solving partial differential equations in which functions are discontinuously approximated within elements. Typically, this is done by using polynomials. Here, the use of radial basis functions as a methods for underlying functional approximation is explored in one dimension, using both analytical and numerical methods. At some mesh densities, RBF flux reconstruction is found to outperform polynomial flux reconstruction, and this range of mesh densities becomes finer as the width of the RBF interpolator is increased. A method which avoids the poor conditioning of flat RBFs is used to test a wide range of basis shapes, and at very small values, the polynomial behaviour is recovered. Changing the location of the solution points is found to have an effect similar to that in polynomial FR, with the Gauss--Legendre points being the most effective. Altering the location of the functional centres is found to have only a very small effect on performance. Similar behaviours are determined for the non-linear Burgers' equation.
\end{abstract}

\begin{keywords}
High Order, Flux Reconstruction, Radial Basis Functions
\end{keywords}


\section{Introduction}\label{sec:intro}
The rapid and accurate numerical solution of partial differential equations (PDEs) is critically important across all disciplines of engineering and the sciences. Key applications include, amongst many others: aerospace engineering~\citep{Temam2001}; astrophyics~\citep{Harris1957}; and geophysical engineering~\citep{Morency2008}. The numerical solution of partial differential equations involves converting the continuous equations into a discrete difference form, which can then be solved approximately. The details of this discrete form --- its derivation --- are what distinguishes between different approaches to solving PDEs. 

In recent years, the attention of the numerical modelling community has been turning towards high order methods, which are thought to offer the potential for computing highly accurate solutions at much reduced computational cost~\citep{Hesthaven2011}. One particularly promising area appears to be the family of Discontinuous Galerkin (DG) or Flux Reconstruction (FR) schemes \citep{Reed1973, Huynh2007}. These schemes have become popular because of their efficiency of implementation on modern computational hardware~\citep{Vermeire2015}, and the elegance of their mathematical representation.

Modern high performance computing facilities have continued to increase their theoretical peak count of floating point operations per second (FLOPS), broadly in line with Moore's Law~\citep{Schaller1997}. However, over the past decade or two these performance gains have largely come not from the increase in individual processor clock speeds, but from the inclusion of increased numbers of computational cores per processor~\citep{top500, Watson2016}, a trajectory which has been taken to extreme levels in clusters based on graphical processing units~\citep{Nvidia2020}. This has caused a new challenge: the memory bandwidth has not kept pace with FLOPS across many architectures, and consequently methods such as FR, which have a high ratio of local computation to remove communication, are well suited to this new computational paradigm~\citep{Iyer2021}.



Alongside the development of these methods, there has been an increased interest in the use of Radial Basis Functions (RBFs) to interpolate functions. Indeed, a number of approaches have used RBFs to solve complex partial differential equations across various disciplines of engineering~\citep{Flyer2007, Larsson2003, Hon1998}. These applications of RBFs generally fall into the category of so-called meshless methods, with two seminal works on the topic introducing the RBF finite difference (RBF-FD) method~\citep{Tolstykh2003} and the meshless Petrov--Galerkin approach~\citep{Atluri1998}. More recently, in the unpublished work of \citet{Glaubitz2021}, a global RBF method that draws upon the family of FR methods has been constructed and proven to be linearly stable. In this paper, we investigate the potential use of RBFs as an underlying method of functional approximation within the broader framework of solving PDEs via Flux Reconstruction.

This investigation is structured as follows. In \cref{sec:rbf}, we introduce RBFs and the techniques used to form both the basis and the subsequent interpolations. \cref{sec:fr} introduces the FR spatial discretisation and detail how RBFs can be applied to FR. In \cref{sec:fourier}, Fourier and combined Fourier analysis are used to investigate the linear stability and effect of RBFs on the numerical properties of FR. In \cref{sec:numerical}, a range of numerical experiments are performed to assess the effects on error and order, which culminate in an investigation of Burgers' turbulence. Finally, the conclusions of this work are drawn in \cref{sec:conclusions}, where insights into the future applications of this methodology are outlined.
\section{Radial Basis Functions}\label{sec:rbf}
    In finite element methods it is common to form an approximation space in some subset of polynomial space. This has several advantages. Polynomials have been widely studied and analysed for centuries, and their strengths and weaknesses are fairly well understood (\textit{e.g.}~\citep{Brutman1978, Gasca2000, Boor1992}). Polynomials are generally quite flexible, and many useful finite element spaces can be formed with them --- such as Lagrange elements, Raviart--Thomas elements, and N\'ed\'elec elements. For most of these methods, it is advantageous to make use of a basis formed from an orthogonal set of polynomials, such as the Jacobi polynomials, $P_n^{\alpha, \beta} \left ( \boldsymbol{x} \right )$. 
    
   One possible alternative to the use of spaces built on polynomials is to make use of RBFs as the basis. RBFs are a class of functions whose value at depends solely on a distance metric --- typically the Euclidean norm --- of the evaluation point, $\boldsymbol{x}$, from the functional centre, $\boldsymbol{r_c}$. This distance metric is usually denoted as a radius, $r = \left \| \boldsymbol{x} - \boldsymbol{r_c} \right \|$. RBFs come in several flavours, and some of the most popular are shown in Table~\ref{tab:rbftypes} below. 
    \begin{table}[hbt!]
        \centering
        \caption{\label{tab:rbftypes}Common Radial Basis Functions.}
        \begin{tabular}{llc}
            \toprule
            RBF & & $\phi(r)$ \\ \midrule \\[-10pt]
            Gaussian & (GA) & $\exp{\left(- \varepsilon^2 r^2\right)}$ \\[5pt]
            Multiquadratic & (MQ) & $\left(1 + \varepsilon^2 r^2\right)^{0.5}$ \\[5pt]
            Inverse Quadratic & (IQ) & $\left(1 + \varepsilon^2 r^2\right)^{-1}$ \\[5pt]
            Inverse Multiquadratic & (IMQ) & $\left(1 + \varepsilon^2 r^2\right)^{-0.5}$ \\[5pt]
            Wendland$_{1,3}$ & (W13) & $\left ( 1 - r \right ) ^4_+ (4r + 1)$ \\[5pt]
            \bottomrule
        \end{tabular}
    \end{table}
    
    Radial Basis Functions have become popular for use in approximating functions both because they are conceptually simple and because they have been shown to be extremely effective at interpolating accurately from randomly scattered points~\citep{Wu1993}. Indeed, there are few theoretical limits to their ability to interpolate over arbitrary sets of points~\citep{Buhmann2003}. Owing to their flexibility, several schemes --- such as a mesh-free Petrov--Galerkin method which made use of the MQ-RBF as a basis~\citep{Atluri1998,Duan2008,Li2019} --- have been developed to exploit their properties. Other works have studied the properties of RBFs in the approximation of functions in various spaces and with reference to Discontinuous Galerkin, such as \citet{Wendland1997} and \citet{Wendland1999}, and these have shown that RBFs can be used to form an effective approximation space under certain conditions.
    
\begin{figure}[htb]
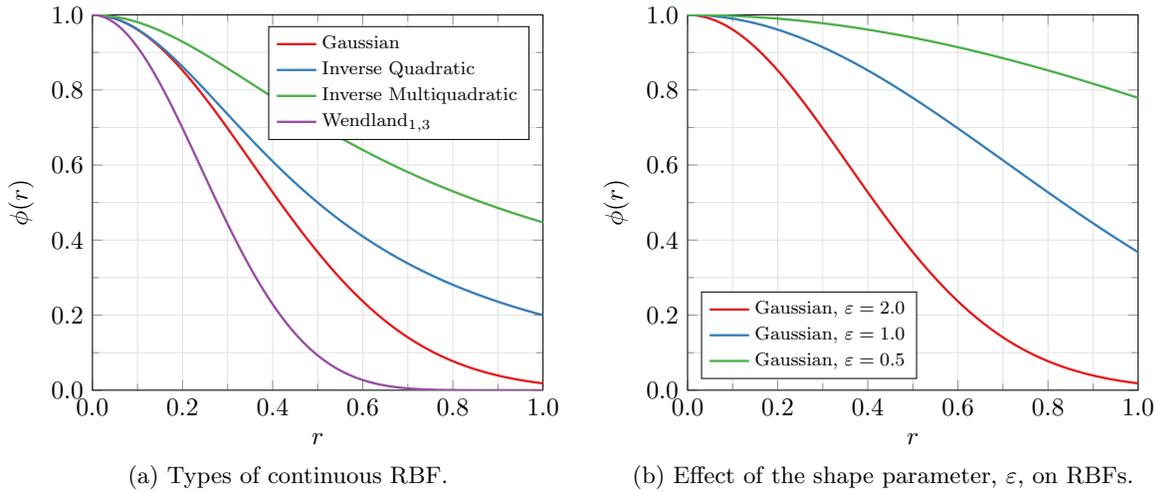

        \centering
        \subfloat[Types of continuous RBF.]{\label{fig:Types_of_RBF}\adjustbox{width=0.49\textwidth,valign=b}{\input{./Figs/tikz/Types_of_RBF_Shape.tex}}}
        ~
        \subfloat[Effect of the shape parameter, $\varepsilon$, on RBFs.]{\label{fig:ShapeFunction_EffectOnRBF}\adjustbox{width=0.49\textwidth,valign=b}{\input{./Figs/tikz/RBF_ShapeParameter_Shape.tex}}}
        \caption{The behaviour of Radial Basis Functions.}
    \end{figure}

    \cref{fig:Types_of_RBF} shows the approximate shape of various commonly seen members of the family of RBFs. All of the functions shown here, with the exception of the Wendland functions~\citep{Wendland1995}, are parameterised by $\varepsilon$, a constant known as the \emph{shape parameter}. Typically, as $\varepsilon$ is reduced, the radial bases become increasingly wide, and tend towards an approximation of $\phi(r) = 1$. The effect of decreasing the shape parameter on the shape of the different radial basis functions is shown for the Gaussian in \cref{fig:ShapeFunction_EffectOnRBF}. 
    
    
    To use radial bases to generate an analytic approximation to a function, $u(\boldsymbol{x})$, which has been evaluated at $N$ points, an $N \times N$ system of equations must be solved:
    
    \begin{align}
        \label{eq:RBF_AlternantMatrix}
        \begin{pmatrix} \phi \left ( \left \| \boldsymbol{x}_1 - \boldsymbol{r_c}_1 \right \| \right) & \phi \left ( \left \| \boldsymbol{x}_1 - \boldsymbol{r_c}_2 \right \| \right) & \dots & \phi \left ( \left \| \boldsymbol{x}_1 - \boldsymbol{r_c}_N \right \| \right) \\ \phi \left ( \left \| \boldsymbol{x}_2 - \boldsymbol{r_c}_1 \right \| \right) & \phi \left ( \left \| \boldsymbol{r_c}_2 - \boldsymbol{x}_2 \right \| \right) & \dots & \phi \left ( \left \| \boldsymbol{x}_2 - \boldsymbol{r_c}_N \right \| \right) \\ \vdots & \vdots & \ddots & \vdots \\ \phi \left ( \left \| \boldsymbol{x}_N - \boldsymbol{r_c}_1 \right \| \right) & \phi \left ( \left \| \boldsymbol{x}_N - \boldsymbol{r_c}_2 \right \| \right) & \dots & \phi \left ( \left \| \boldsymbol{x}_N - \boldsymbol{r_c}_N \right \| \right) \end{pmatrix} \begin{pmatrix} \boldsymbol{\psi}_1 \\ \boldsymbol{\psi}_2 \\ \vdots \\ \boldsymbol{\psi}_N \end{pmatrix} = \begin{pmatrix} u(\boldsymbol{x}_1) \\ u(\boldsymbol{x}_2) \\ \vdots \\ u(\boldsymbol{x}_N) \end{pmatrix}
    \end{align}
    which is analogous to using a Vandermonde matrix to compute polynomial functional approximations, and which results in the functional approximation:
    
    \begin{align}
    \label{eq:RBF_Decomposition}
    u(\boldsymbol{x}) \approx \sum^{N}_{i=1} \psi_i \phi \left ( \left \| \boldsymbol{x} - \boldsymbol{r_c}_i \right \| \right) 
    \end{align}
    
    As the shape parameter is decreased towards zero, the alternant matrix governing the system of equations in \ref{eq:RBF_AlternantMatrix} approaches the matrix of ones, and so becomes increasingly ill-conditioned. However, there are several ways to avoid this ill-conditioning, which will be discussed more extensively later in the paper.

\section{Flux Reconstruction}\label{sec:fr}
To provide some necessary background, the basic procedure for Flux Reconstruction as applied to the linear advection--diffusion equation, as originally developed by \citet{Huynh2007}, is outlined in this section. The approach taken here largely follows that outlined by \citet{Witherden2014}. To begin, the linear advection--diffusion equation is written in strong conservation form:
\begin{align}
    \label{eq:advectionDiffusionEquation}
    \px{\boldsymbol{u}}{t} + \px{\boldsymbol{f}_i}{x_i} = 0
\end{align}
where:
\begin{align}
    \boldsymbol{u} = \begin{pmatrix} u \end{pmatrix} \quad , \quad \boldsymbol{f} = \begin{pmatrix} a_i u - \mu_i \frac{\partial u}{\partial x_i} \end{pmatrix}
\end{align}

The computational domain, $\Omega$, is partitioned into $N$ compatible non-overlapping elements, $\Omega_n$, such that:
\begin{align}
    \Omega = \bigcup^N_{n=1} \Omega_n \qquad \text{and} \qquad \bigcap^N_{n=1} \Omega_n = \emptyset
\end{align}

\begin{figure}[htbp]
    \centering
        \adjustbox{width=0.6\linewidth, valign=b}{\input{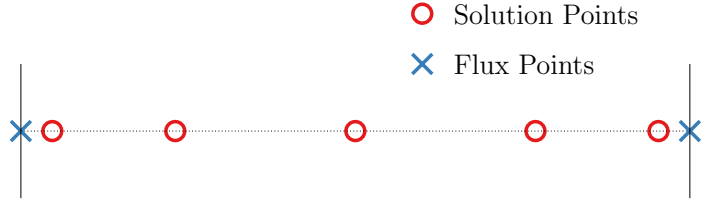}}
    \caption{\label{fig:SolutionAndFluxLayout}A one dimensional example of solution point and flux point layouts, $n_s=5$}
\end{figure}

Each of these subdomains is, for computational efficiency, mapped onto a standard reference domain, $\Omega_n \to \Omega_s$, in which $\Omega_s \in \mathbb{R}^d$, where $d$ denotes the dimensionality of the problem. For example, the typical computational reference domain used for quadrilateral elements is $\Omega_s=[-1,1]^2$. Solution points are placed at specified locations within the domain, and flux points around its edges, as shown in \cref{fig:SolutionAndFluxLayout} for a one dimensional element. To formalise the transform between the \emph{real} and \emph{transformed} domains, $\Omega_n \to \Omega_s$, a Jacobian matrix, $\boldsymbol{J}_n$, is defined:
\begin{subequations}
    \begin{align}
        \boldsymbol{u} = \boldsymbol{u}(\hat{\boldsymbol{x}}, t) = \left | \boldsymbol{J}_n \right | u^\delta (\boldsymbol{x}, t) \\
        \hat{\boldsymbol{f}}^\delta = \hat{\boldsymbol{f}}^\delta (\hat{\boldsymbol{x}}, t) = \left | \boldsymbol{J}_n \right | \boldsymbol{J}_n^{-1} \boldsymbol{f}^\delta (\boldsymbol{x}, t)
    \end{align}
\end{subequations}
where the hat indicates transformed variables.

The problem is assumed to be formulated as an initial value problem, so that the initial solution at $t=0$ is known at all the solution points. The conserved variables are approximated within each subdomain as a weighted sum of basis functions:
\begin{align}
    \label{eq:SumofBasisFunctions}
    \hat{u}^\delta (\hat{\boldsymbol{x}}) = \sum^{n_s}_{i=1} \hat{u}_i^\delta \theta_i(\hat{\boldsymbol{x}})
\end{align}

In one spatial dimension a Lagrange polynomial, for example, can be used as the basis function:
\begin{align}
    \label{eq:LagrangeInterpolation}
    \theta_k (\hat{x}) = \prod^{n_s}_{i=1,i \ne k} \frac{\hat{x} - \hat{x}_i}{\hat{x}_k - \hat{x}_i}
\end{align}

When working in higher dimensions, it is straightforward to extend this approach to tensor-product elements with maximal order bases, but other bases are better suited for simplices, such as the Dubiner basis for triangles~\citep{Dubiner1991}.

Once the conserved variables are functionally approximated within the subdomain based on their solution point values, these bases are projected to the flux points. This allows the values of the conserved variables at the flux points to be determined:
\begin{align}
    \label{eq:SolutionToFluxConserved}
    \hat{u}^\delta_F = \hat{u}^\delta (\hat{\boldsymbol{x}}_F)
\end{align}

The flux points of the compatible subdomains are colocated with those of their neighbours. As the projection from the solution points has been performed separately for each subdomain, the two sets of values of the conserved variables at each pair of colocated flux points are likely to be different.

The gradient of the conserved variables at each of the solution points is calculated in the transformed domain by simply differentiating the bases:
\begin{align}
    \px{\hat{u}^\delta}{\hat{x}_j} = \sum^{n_s}_{i=1} \hat{u}_i^\delta \px{\theta_i(\hat{x})}{\hat{x}_j}
\end{align}

This gradient is formed from a continuous solution within each element, however, in order for the resulting scheme to be conservative between elements, it is necessary to correct the gradient using a common interface solution, $\hat{u}^c$, and then to propagate this into both cells by using a correction function, $g$. In one dimension, incorporating the corrections from the colocated flux points at both the left and right sides of the cell, the process is given by:
\begin{align}
    \px{\hat{u}}{\hat{x}} = \sum^{n_s}_{i=1} \hat{u}_i^\delta \px{\theta_i(\hat{x})}{\hat{x}} + \left ( \hat{u}^c_L - \hat{u}^\delta_{F,L}\right ) \dx{g_L(\hat{x})}{\hat{x}} + \left ( \hat{u}^c_R - \hat{u}^\delta_{F,R}\right ) \dx{g_R(\hat{x})}{\hat{x}}
\end{align}

Applying this correction ensures that gradient approximation is in the point-wise $H^0_\mathrm{div}$ space, which is a regularity constraint that can be derived from the weak form of the equations.

A number of approaches for calculating the common interface values exist~\citep{Huynh2009}, but here a simple arithmetic mean of the two discontinuous flux point values at the collocated interface is suggested, so that, in the $i^{th}$ cell:
\begin{align}
    \hat{u}^c_{L,i} = \frac{1}{2} \left ( \hat{u}^\delta_{F,L,i} + \hat{u}^\delta_{F,R,i-1} \right ) \qquad \text{and} \qquad \hat{u}^c_{R,i} = \frac{1}{2} \left ( \hat{u}^\delta_{F,R,i} + \hat{u}^\delta_{F,L,i+1} \right )
\end{align}

Now, the flux at each of the solution points can be calculated in the computational domain from the conserved variables, $\hat{\boldsymbol{f}}^\delta_s = \hat{\boldsymbol{f}}^\delta \left ( \hat{u}_s^\delta, \frac{\partial \hat{u}}{\partial \hat{x}_j} \big\vert_s \right )$, with the advective-diffusive fluxes being dependent on both the conserved variables themselves and their corrected gradients.

Once there is a functional approximation of the fluxes within the subdomain, these can be evaluated at the flux points, and simultaneously projected in the direction of the computational normal at the corresponding flux point to get the discontinuous normal computational fluxes, $\left ( \hat{\boldsymbol{f}} \cdot \hat{\boldsymbol{n}} \right )^\delta_F = \hat{F}^\delta_n $.

To ensure conservation is enforced between cells, a common interface normal flux is calculated at the collocated flux points, $\left ( \hat{\boldsymbol{f}} \cdot \hat{\boldsymbol{n}}\right )^c_F = \hat{F}^c_n$. As with the common interface value, there is a choice in how these fluxes are computed. Typically, the fluxes are divided into advective and diffusive components, such that $\left ( \hat{\boldsymbol{f}} \cdot \hat{\boldsymbol{n}}\right )^c_F = \hat{F}^c_n = \hat{F}^c_{n,\text{Adv}} + \hat{F}^c_{n,\text{Dif}}$. The contribution of the advective fluxes can then be estimated with an approximate Riemann solver, such as a Roe flux~\citep{Roe1981}, while a typical approach for the diffusive is taking the arithmetic mean of those calculated on the two discontinuous sides. When this approach is coupled to the use of the arithmetic mean for the common interface value, it is generally referred to as BR2~\citep{Bassi2000}.

Equipped with the common interface normal fluxes, the divergence of the fluxes can estimated within each cell in a manner which enforces conservation between cells. This is given in one dimension by:
\begin{align}
    \label{eq:FR_advectionDiffusionEquation}
    \frac{\partial \hat{\boldsymbol{f}}}{\partial \hat{x}} = \sum^{n_s}_{i=1} \hat{\boldsymbol{f}}_i^\delta \frac{\partial \theta_i(\hat{x})}{\partial \hat{x}} + \left ( \hat{F}^c_{n,L} - \hat{F}^\delta_{n,L}\right ) \frac{\partial g_L(\hat{x})}{\partial \hat{x}} + \left ( \hat{F}^c_{n,R} - \hat{F}^\delta_{n,R}\right ) \frac{\partial g_R(\hat{x})}{\partial \hat{x}}
\end{align}

With the approximation of the divergence of the fluxes from~\cref{eq:advectionDiffusionEquation} given by~\cref{eq:FR_advectionDiffusionEquation}, it is possible to advance the conserved variables in time through the use of a classical ODE integrator, such as a Runge--Kutta scheme.

\subsection{Correction Functions for FR}\label{ssec:correctionfunctions}
    One of the defining features of flux reconstruction is its use of a correction function, $g$, to tie the individual cells together, thereby enforcing intercell continuity. The shape of this correction function has been shown to have a dramatic influence on the behaviour of the scheme~\citep{Vincent2011}. 
    
    There have been several investigations into novel families of correction functions, usually based on bounding the growth of some kind of energy norm~\citep{Castonguay2013,Vincent2010,Trojak2019}. The resulting schemes have proven to be highly effective in improving accuracy and maximising stability, but can be tricky to extend to non-tensor product shapes, because of their complex structures in higher dimensions. For example, the VCJH scheme~\citep{Vincent2011} on simplices in two dimensions:
    \begin{equation}
        \label{eq:VCJHcorrection}
        c \sum^{n_s}_{k=1} \sigma_{f,k} \sum^{p+1}_{m=1} \begin{pmatrix} p \\ m - 1 \end{pmatrix}  \left ( D^{(m,p)} \xi_i \right ) \left ( D^{(m,p)} \xi_k \right ) = -\sigma_{f,i} + \int_{\boldsymbol{\Gamma}_S} \left ( \boldsymbol{g}_f \cdot \bar{\hat{\boldsymbol{n}}} \right ) \xi_i d\Gamma
    \end{equation}
    where the operator $D^{(m,p)}$ is defined as:
    \begin{equation*}
        D^{(m,p)} = \frac{\partial^p}{\partial r^{p-m+1}\partial s^{m-1}}
    \end{equation*}
    and where $r$ and $s$ are some orthogonal coordinates. The coefficients of the divergence of the correction function corresponding to the $f^\mathrm{th}$ flux point are given by \cref{eq:VCJHcorrection}, such that:
    \begin{equation}
        \label{eq:correctionFunctionSum}
        \nabla \cdot \boldsymbol{g}_f = \sum^{n_s}_{j=1} \sigma_{f,j} \xi_j \left ( \hat{\boldsymbol{x}} \right )
    \end{equation}
    where $\xi_j$ forms an orthonormal basis with support on the reference domain.
    
    When the parametrisation constant, $c$, is set to zero, this function recovers a colocated nodal Discontinuous Galerkin scheme, but non-zero values can give enhanced performance. It has also been shown that modifying the correction function can be thought of as applying a filter to the collocated nodal Discontinuous Galerkin scheme~\citep{Zwanenburg2016}. However, when setting $c$ to zero in \cref{eq:VCJHcorrection} for computing the VCJH correction function, it collapses to the far simpler:
    \begin{equation}
        \label{eq:VCJHcorrection_DG}
        \sigma_{f,i} = \int_{\boldsymbol{\Gamma}_S} \left ( \boldsymbol{g}_f \cdot \bar{\hat{\boldsymbol{n}}} \right ) \xi_i d\Gamma
    \end{equation}
    which gives the required coefficients for \cref{eq:correctionFunctionSum} directly. A further step eliminates the dependency of the right hand side on the correction function by noting that, to ensure continuity, $\left ( \boldsymbol{g}_f \cdot \bar{\hat{\boldsymbol{n}}} \right )$ must equal unity at flux point $f$ and zero at the others. This allows the term to be replaced by a flux point identity basis, which takes a value of one at its own flux point, and zero at every other:
    \begin{equation}
        \zeta^F_k \left ( \hat{\boldsymbol{x}}_f \right ) = \delta_{f,k}
    \end{equation}
    
    This finally results in an easily solvable set of coefficients for the correction functions:
    \begin{equation}
        \label{eq:VCJHcorrection_DG_simp}
        \sigma_{f,i} = \int_{\boldsymbol{\Gamma}_S} \zeta^F_f \xi_i d\Gamma
    \end{equation}
    which simply require knowledge of an orthonormal basis over the transformed domain, and an identity basis over the flux points.
    
    No requirement is placed on the functional form of the bases, $\xi_j$, used to compute the correction function field --- except for the requirement for orthonormality. However, the Gram--Schmidt process provides a tool to generate an orthonormal set of bases from any spanning set.

\subsection{Summation-by-parts}\label{ssec:sbp}

Before introducing the application of RBFs to FR, it is worth first taking a brief diversion into some approaches for testing the stability of the numerical method. From the work of \citet{Ranocha2016} and others, the FR method can be cast into the summation-by-parts (SBP)  framework~\citep{Strand1994,DelReyFernndez2014}. Here, a one dimensional definition of the SBP framework is used, although higher dimensionality versions are also available~\citep{DelReyFernndez2014,Hicken2016}.
    
\begin{definition}[Summation-by-parts]
    For a set of operators, $\mathbf{M}, \mathbf{D}, \mathbf{P}$, and  $\mathbf{B}$, defined by:
    \begin{subequations}\label{eq:sbp}
        \begin{align}
            \int_{\hat{\Omega}} \hat{v}\hat{w}\mathrm{d}x &\approx \langle\hb{v},\hb{w}\rangle_M \equiv \hb{v}^T\mathbf{M}\hb{w} \\ \px{\hat{w}}{x}\bigg|_{x=\hb{x}} &\approx \mathbf{D}\hb{w} \\
            \mathbf{B} &= \mathrm{diag}{(-1,1)} \\
            \mathbf{P}\hb{w} &\approx [\hb{w}_L, \hb{w}_R]^T
        \end{align}
    \end{subequations}
    then these operators are said to be summation-by-parts operators if:
    \begin{subequations}
        \begin{align}
            \int_{\hat{\Omega}}w\px{v}{x}\mathrm{d}x + \int_{\hat{\Omega}} v\px{w}{x}\mathrm{d}x = vw\big|_{\partial\hat{\Omega}} &\approx \hb{w}^T\mathbf{MD}\hb{v} + \hb{w}^T\mathbf{D}^T\mathbf{M}\hb{v}\\
            \mathbf{MD} + \mathbf{D}^T\mathbf{M} &= \mathbf{P}^T\mathbf{BP}
        \end{align}
    \end{subequations}
\end{definition}
        
From \citet{Ranocha2018}, we state two lemmas on the conservation and stability of the FR method for linear advection equations:
\begin{lemma}[Conservation]\label{thm:sbp_cons}
    Given a set of SBP operators, $\{\mathbf{M}, \mathbf{D}, \mathbf{P}, \mathbf{B}\}$ and the FR lifting operator, $\mathbf{C}$, then for the one dimensional linear advection equation this set is conservative if:
    \begin{equation}
        \mathbf{1}^T\mathbf{MC} = \mathbf{1}^T\mathbf{P}^T\mathbf{B} 
    \end{equation}
\end{lemma}
        
\begin{lemma}[Linear Stability]\label{thm:sbp_stab}
    Given a set of SBP operators, $\{\mathbf{M}, \mathbf{D}, \mathbf{P}, \mathbf{B}\}$, and the FR lifting operator, $\mathbf{C}$, then for the one dimensional linear advection equation this set is linearly stable if:
    \begin{equation}
        \mathbf{C} = \mathbf{M}^{-1}\mathbf{P}^T\mathbf{B}
    \end{equation}
\end{lemma}

These lemmas prove to be useful in testing the linear stability and conservation of schemes by casting the operators in the form of \cref{eq:sbp}, and will be used later.
    
\subsection{RBFs for Flux Reconstruction}\label{ssec:method}

    RBFs can be used to form the basis of the approximation within the transformed cell in the same way as polynomials. For computational efficiency, most of the functional approximation and evaluation operations in FR are done on the transformed elements, $\boldsymbol{\Omega}_s$, where they can be expressed as matrix operations of the form:
    
    \begin{align}
        \hat{u}^\delta_F = \mathbf{M}_1 \hat{u}^\delta_S
    \end{align}
    which applies the $\mathbf{M}_1$ operator to compute the flux point conserved variables from the solution point conserved variables --- steps~\cref{eq:SumofBasisFunctions,eq:LagrangeInterpolation,eq:SolutionToFluxConserved} above. Critically, the operators themselves do not change during the calculation. To find $\mathbf{M}_1$, systems of equations are solved \emph{a priori} on the transformed domain, $\Omega_s$, and the required matrices are generated, which can then be applied at each time step.
    
    Making the change from FR based on polynomial approximation to FR based on radial basis approximation essentially involves swapping out the polynomial operators and replacing them with operators derived using radial bases. The RBF-derived operators themselves are computed by manipulating the alternant matrix of~\cref{eq:RBF_AlternantMatrix}.
    
    \begin{figure}[htb]
        \centering
        \adjustbox{width=0.8\linewidth, valign=b}{\input{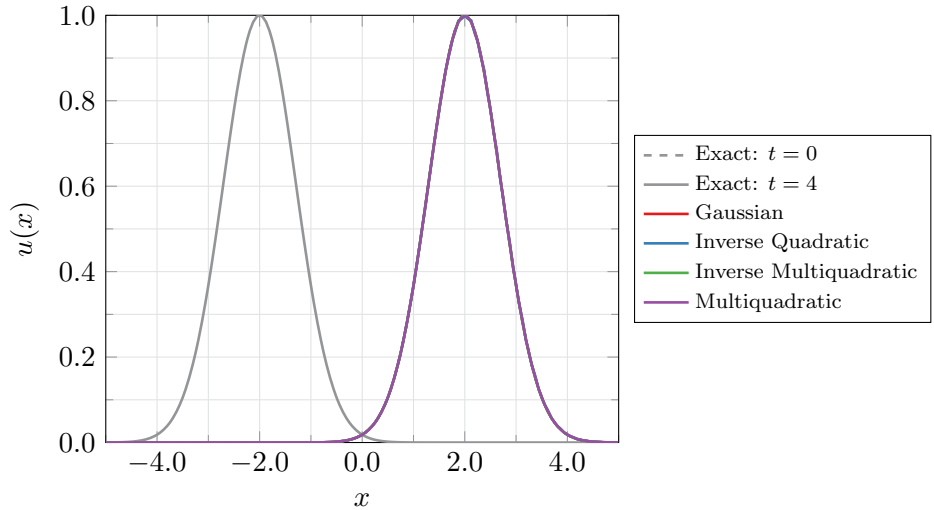}}
        \caption{\label{fig:AllRBF_AdvGaussian}Advection of a Gaussian with RBF functional approximation.}
    \end{figure}
    
    \cref{fig:AllRBF_AdvGaussian} shows the use of various different radial bases, along with the baseline polynomial case, to advect and diffuse a Gaussian bump on a one dimensional domain. It can be seen in the figure that all of the different bases used are able to capture the basic behaviour of the wave as it travels across the domain. This suggests that RBFs are in principle able to serve as the underlying functional approximation method for flux reconstruction. Over the following sections, the behaviour and performance of radial basis function flux reconstruction (RBF-FR) will be explored in more detail.
    
    To keep the scope manageable, the rest of this paper will focus solely on the use of Gaussian RBFs, where the interpolating function is given by $\phi(r) = \exp{(-\varepsilon^2 r^2)}$. The Gaussian is chosen because of its familiarity, its ease of manipulation, and because it has been shown to have some particularly intriguing properties in the flat limit~\citep{Fornberg2004a}, as the shape parameter $\varepsilon \to 0$.
  
\subsection{Investigating RBF-FR}

By changing the method of approximation within the transformed domains from polynomials to RBFs, various modifications to the simulation are possible. Here, we detail the modifications investigated in this work.

\begin{description}
\item[Solution Point Locations.]
    The layout of the solution points within the reference cell is known to have an influence on the performance of flux reconstruction~\citep{Witherden2014a, Witherden2021}. 

    \begin{figure}[htbp]
        \centering
            \adjustbox{width=0.6\linewidth, valign=b}{\input{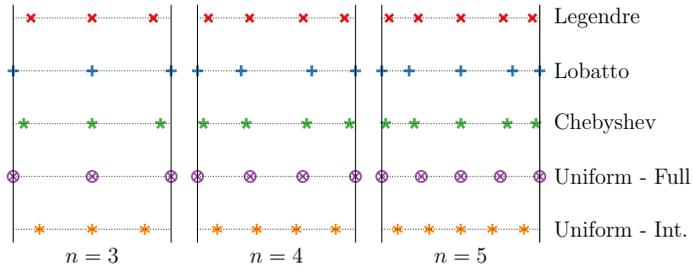}}
        \caption{\label{fig:pointLayouts}Canonical point layouts within the reference element.}
    \end{figure}

    \cref{fig:pointLayouts} shows the one dimensional solution point locations that have been tested here. The Legendre, Lobatto, and Chebyshev points are familiar from polynomial quadrature theory as the roots of the relevant functions. They have different advantages and disadvantages in terms of function approximation, and the Lobatto points in particular have some computational advantages in terms of FR. As seen in the figure, the Lobatto points include the end points. This means that the flux points are colocated with the terminal solution points, avoiding the need for projection to the flux points, and thereby offering a saving in compute time. In addition to these three familiar point layouts, two others have been included --- one, the `internal' uniform points, which are uniformly distributed over the inside of the cell, with half a spacing between the terminal solution points and the flux points, and two, the `full' uniform points, which are uniformly distributed over the entire cell, from flux point to flux point, and which also share the potential computational benefit of colocated flux and solution points with the Lobatto points.
    
    Although distinct points are not able to form non-unisolvent sets in one dimension, this is something which may occur in higher dimensions. RBFs are more geometrically flexible than polynomials, so in two and higher dimensions, they are able to interpolate a wider number of scattered points. Two and higher dimensional cases are outside the scope of the present paper, but form a strong motivation for the development of RBF-FR as outlined here.

\item[RBF Functional Centre Locations.]
    Along with their shape, radial basis functions are defined by the locations from which the distance $r$ is measured --- the functional centres of the bases. Frequently in the use of RBFs, the functional centres are chosen to be colocated with the measuring points, but this is not a requirement. To investigate the influence of different sets of functional centres, the same set of points shown in \cref{fig:pointLayouts} are used.

\item[Shape Parameter.]
    The shape parameter, $\varepsilon$, or equivalently, the width of the basis, is known to have a significant effect on the performance of RBF interpolation. Indeed, a number of investigations have been undertaken to tune the value of the shape parameter to give the best interpolating performance~\citep{Fasshauer2007,Bayona2011,Davydov2011}, and so the influence of this parameter is explored here. It is also possible to have different bases take different shape parameters within the same functional interpolation. This kind of local variation in $\varepsilon$ is not considered further within this paper, but it remains something which may be of interest for future research.

\item[Small Values of Shape Parameter.]
    Although the value of the shape parameter is critical to the performance of an RBF interpolation, the RBF alternant matrix shown in \cref{eq:RBF_AlternantMatrix} becomes increasingly ill-conditioned as the the shape parameter decreases, becoming equal to the matrix of ones in the limit $\varepsilon \to 0$. \cref{fig:condNumber_RBFDirect} shows this variation in condition number with shape parameter. The condition number is also more weakly affected by the location of the functional centres, an effect which is related to the mean distance between them. A general rule-of-thumb is that a for every order of magnitude the condition number increases by, a digit of precision is lost in the solution. This means that the system of equations cannot reliably be solved for smaller values of shape parameter. However, it is not the basis \emph{itself} that is ill-conditioned --- rather, it is an intermediate step in calculating it which is to blame.

    \begin{figure}[tbhp]
        \centering
        \adjustbox{width=0.6\linewidth, valign=b}{\input{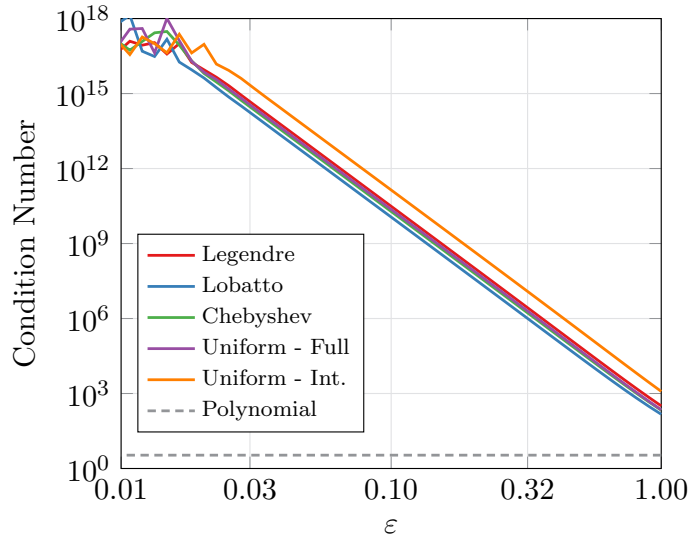}}
        \caption{\label{fig:condNumber_RBFDirect} Effect of reducing $\varepsilon$ on alternant matrix condition number when using direct Gaussian RBFs.}
    \end{figure}

    To this end, various methods have been developed which circumvent this step to allow an effective basis to be produced when using very flat radial bases. Examples of ways in which this has been achieved are the Contour-Padé method~\citep{Fornberg2004} and the RBF-QR method~\citep{Fornberg2011}. Here, since the investigation is restricted to Gaussian radial bases, an approach known as the RBF-GA method~\citep{Fornberg2013} is used to build the basis functions without the need to directly solve the system of equations. The details of the approach are omitted here for brevity, but it involves expanding the original Gaussian functions into a Taylor series, and then selecting linear combinations which cause most of these terms to cancel out. The result is a new basis which, critically, spans exactly the same space as the original, but which is well conditioned even at very low values of $\varepsilon$.

    \cref{fig:condNumber_RBFGA} shows the effect on the condition number of the solution point alternant matrix of calculating the basis via the RBF-GA method, rather than the RBF-Direct approach. The same dependence on the mean spacing between the functional centres is observed, but in general the condition number of the alternant is below $10^3$, meaning that by side-stepping the direct solution, it becomes numerically feasible to use far smaller values of the shape parameter.

    \begin{figure}[tbhp]
        \centering
        \adjustbox{width=0.80\linewidth, valign=b}{\input{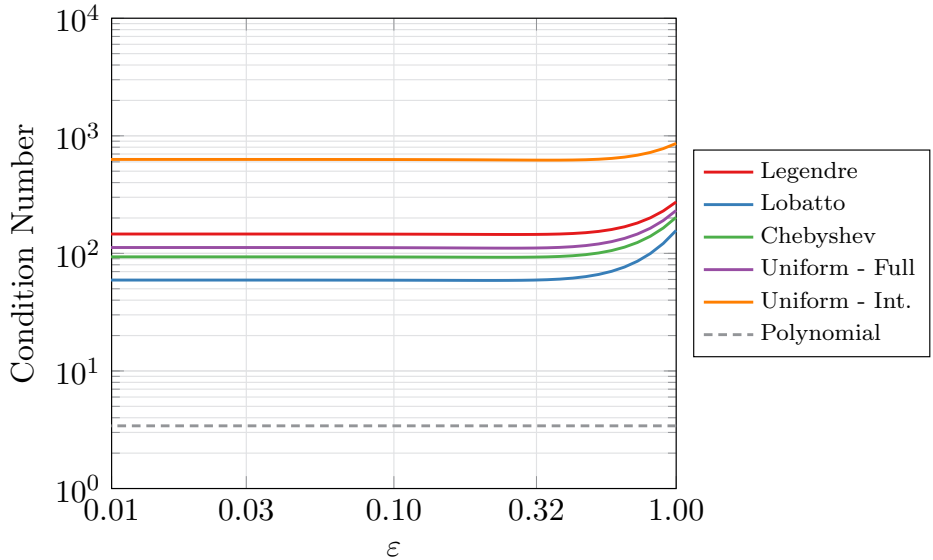}}
        \caption{\label{fig:condNumber_RBFGA} Effect of reducing $\varepsilon$ on alternant matrix condition number when using RBF-GA.}
    \end{figure}
\end{description}
\section{\label{sec:fourier}Fourier Analysis}
    Initially, the properties of RBF-FR schemes are examined analytically. Subsequent sections will illustrate their behaviour through numerical experiments. Fourier analysis provides a powerful tool for examining spatial operators through their frequency responses. Previously in the FR literature, Fourier analysis has been used to understand the dispersion and dissipation of schemes in one and two dimensions for advection equations~\citep{Huynh2007,Trojak2020}. Further works have explored schemes' CFL limits, diffusion effects, and rates of convergence~\citep{Vincent2011,Watkins2016,Asthana2017}. 
    
    As discussed, the FR method generates a solution through the superposition of several modes. When viewed in discrete Fourier space, these have often been divided into primary and secondary modes or into physical and spurious modes. The primary mode describes the most energised, whereas the physical mode is that which most accurately represents the projection from the approximation space into Fourier space. This approach has provided considerable insight into the behaviour of spatial methods. However, an alternative `holistic' approach, in which the combined effect of all the modes is considered, has recently been presented \citep{Alhawwary2020,Asada2020,Fernandez2019}.
    
    We begin by introducing the discrete eigenmode approach. Take the linear advection equation:
    \begin{equation}\label{eq:lin_adv}
        \px{u}{t} + \px{u}{x} = 0, \quad \forall u(x,t): \Omega\times\mathbb{R}_+ \mapsto \mathbb{C}
    \end{equation}
    
    Here, we allow the solution to be complex. If this is taken as a Cauchy problem --- by setting the domain $\Omega$ to be either toroidal or infinite --- then \cref{eq:lin_adv} can support Bloch wave solutions:
    \begin{subequations}
        \label{eq:bloch}
        \begin{align}
            u &= \exp{\left(ik(x - ct)\right)} \\
            \boldsymbol{u}_j &= \exp{\left(ik\left(x_j + \frac{h}{2}(\pmb{\xi}+1) + ct\right) \right)} \quad \forall \boldsymbol{u}_j : \boldsymbol{x}_j \times \mathbb{R}_+ \mapsto \mathbb{C}^n
        \end{align}
    \end{subequations}
    where $\boldsymbol{x}_j \in \Omega^n_j$ for $\Omega_j = [x_j,x_{j+1}]$, where the mesh is taken to be uniform with spacing $h=x_{j+1}-x_j$. The latter equation here represents the discrete solution vector within the $j^\mathrm{th}$ element.
    
    Constructing the flux reconstruction operator, \cref{eq:lin_adv} is first written as:
    \begin{equation}\label{eq:linadv_approx}
        \px{\boldsymbol{u}_j}{t} = -\px{\boldsymbol{u}_j}{x} \approx \mathbf{Q}\boldsymbol{u}_j
    \end{equation}
    then, following the method of \citet{Trojak2020}, we have:
    \begin{equation}
            \mathbf{Q} = -\frac{2}{h}\left(\mathbf{C}_-\exp{(-ikh)} + \mathbf{C}_0 + \exp{(ikh)}\mathbf{C}_+\right)
    \end{equation}
    where:
    \begin{subequations}
        \begin{align}
            \mathbf{C}_{-} = \mathbf{C}\mathbf{K}_{-}\mathbf{P}^T, \quad &\mathbf{K}_{-} = \begin{bmatrix}0 & 0 \\ 1-\alpha & 0 \end{bmatrix} \\
            \mathbf{C}_{0} = \mathbf{D} + \mathbf{C}\mathbf{K}_{0}\mathbf{P}^T, \quad &\mathbf{K}_{0} = -\diag{(\alpha,1-\alpha)} \\
            \mathbf{C}_{+} = \mathbf{C}\mathbf{K}_{+}\mathbf{P}^T, \quad &\mathbf{K}_{+} = \begin{bmatrix}0 & \alpha \\ 0 & 0 \end{bmatrix}
       \end{align}
    \end{subequations}
    
    Here, $\alpha$ is the interface upwinding parameter, with $\alpha=1$ providing a fully upwinded interface and $\alpha=0.5$ giving centrally differenced interfaces.
    
    By substituting \cref{eq:bloch} into \cref{eq:linadv_approx} and rearranging, the modified angular frequency is obtained as:
    \begin{equation}
        c^\prime k\mathbf{w} = \omega^\prime\mathbf{w} = i\mathbf{Q}\mathbf{w}
    \end{equation}
    where $\mathbf{w}$ is an eigenvector, and $c^\prime$ and $\omega^\prime$ are the modified wavespeed and modified angular frequency, respectively. The implication of this is that for a given wavenumber, $k$, the solution to the Bloch wave problem actually given by an FR scheme is equivalent to:
    \begin{equation}
        u_\text{mod}=\exp{(ikx-i\omega^\prime t)}
    \end{equation}

\begin{figure}[tbhp]
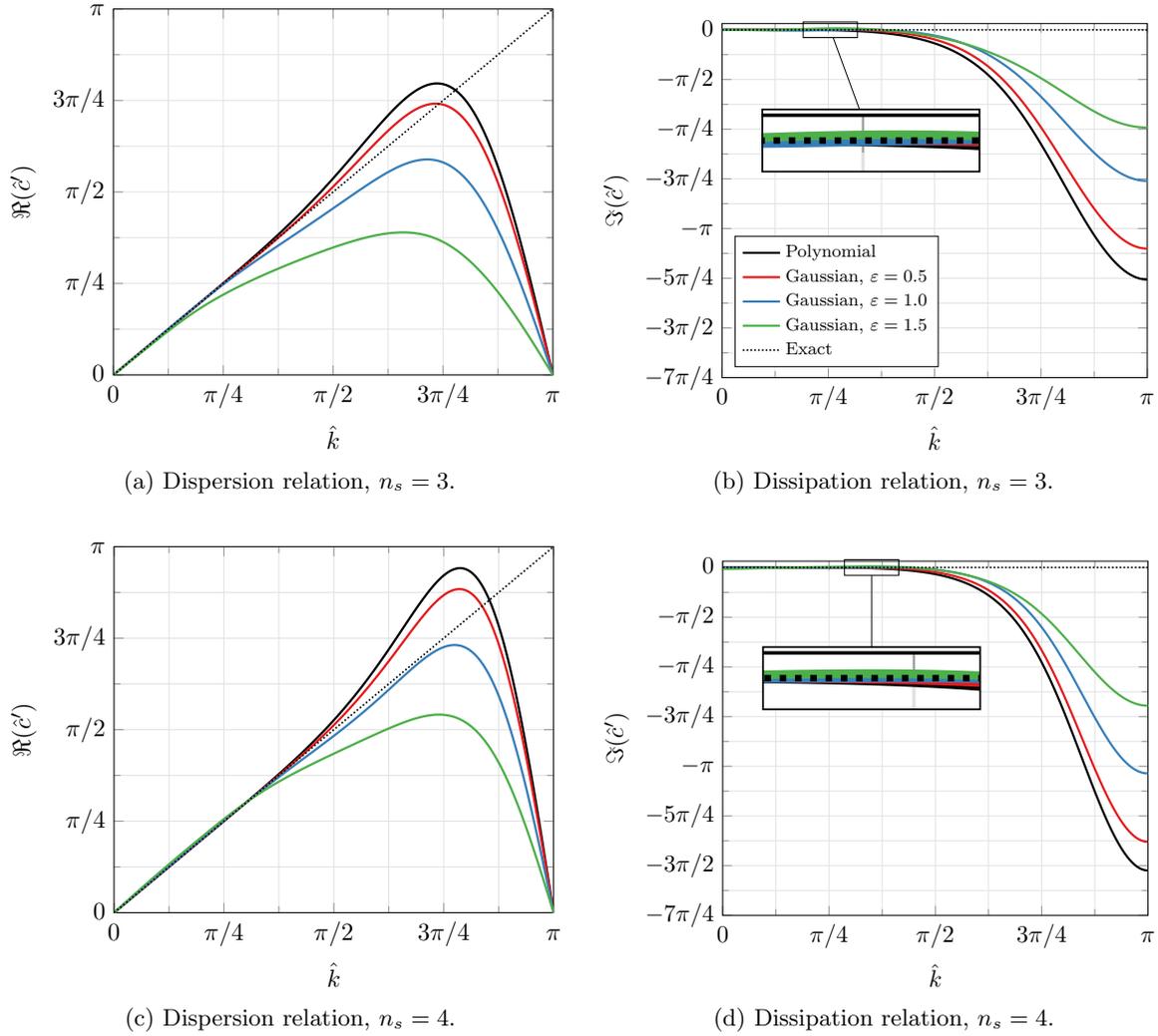

    \centering
    \subfloat[Dispersion relation, $n_s=3$.]{\label{fig:n3_disp}\adjustbox{width=0.49\linewidth,valign=b}{\input{Figs/tikz/n3_gauss_dispersion}}}
    ~
    \subfloat[Dissipation relation, $n_s=3$.]{\label{fig:n3_diss}\adjustbox{width=0.49\linewidth,valign=b}{\input{Figs/tikz/n3_gauss_dissipation}}}
    \newline
    \subfloat[Dispersion relation, $n_s=4$.]{\label{fig:n4_disp}\adjustbox{width=0.49\linewidth,valign=b}{\input{Figs/tikz/n4_gauss_dispersion}}}
    ~
    \subfloat[Dissipation relation, $n_s=4$.]{\label{fig:n4_diss}\adjustbox{width=0.49\linewidth,valign=b}{\input{Figs/tikz/n4_gauss_dissipation}}}
    \caption{\label{fig:disp_diss}Upwind RBF-FR dispersion and dissipation plots.}
\end{figure}

The results presented in \cref{fig:disp_diss} show the real and imaginary components of this modified angular frequency, which correspond to the dispersion and dissipation, respectively, for an RBF-FR scheme. For these calculations, the solution points have been placed at the Gauss--Legendre points, as have the centres of the radial bases. Only the physical mode is shown. These results show that the change of basis has a significant impact on the dispersion and dissipation characteristics, which trend towards the behaviour of the polynomial basis scheme as the shape parameter, $\varepsilon$, is reduced and the RBFs flattened. Furthermore, as $\varepsilon$ is increased the characteristic dispersion over shoot and high dissipation at high wavenumbers is reduced. Which is most likely linked to a flattening of the correction function as $\varepsilon$ increases. Ultimately, continued increases to the shape parameter lead the imaginary part to become positive for some wavenumbers, suggesting a degree of instability.

\begin{figure}[tbhp]
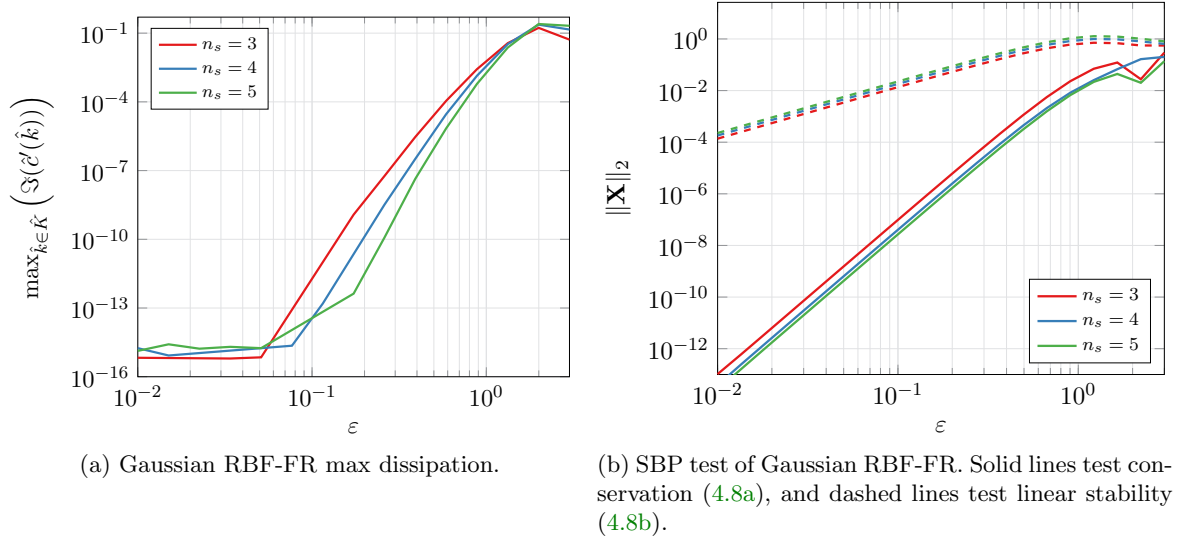

    \centering
    \subfloat[Gaussian RBF-FR max dissipation.]{\label{fig:diss_max}\adjustbox{width=0.49\linewidth, valign=b}{\input{Figs/tikz/dmax_plot}}}
    ~
    \subfloat[SBP test of Gaussian RBF-FR. Solid lines test conservation (\ref{eq: SBP_conservation}), and dashed lines test linear stability (\ref{eq: SBP_stability}).]{\label{fig:sbp}\adjustbox{width=0.49\linewidth, valign=b}{\input{Figs/tikz/sbp_eps}}}
    \caption{Testing the stability of Gaussian RBF-FR schemes.}
\end{figure}

Since changing the RBF width can lead to significant changes in the numerical properties of the method, a sweep was taken over a range of $\varepsilon$ values, and the maximum value of the dissipation recorded. The results for Legendre functional centres and Legendre solution points with a Gaussian RBF are presented in \cref{fig:diss_max}. This shows that the change in stability with $\varepsilon$ is dependent on degree of the basis, with higher degree bases permitting higher $\varepsilon$ values.

Pursuing this further, the SBP framework introduced in \cref{ssec:sbp} gives a means to test for conservation and linear stability in a manner which is independent of the initial condition used. Using this, the norm of the error has been calculated and plotted in \cref{fig:sbp} for Legendre solution points and centres with Gaussian RBFs of varying width. The two norms investigated are:
\begin{subequations}
    \begin{align}
        \|\mathbf{1}^T\mathbf{MC}-\mathbf{1}^T\mathbf{P}^T\mathbf{B}\|_2 &= \text{Conservation Error} \label{eq: SBP_conservation}\\
        \|\mathbf{M}^{-1}\mathbf{P}^T\mathbf{B}-\mathbf{C}\|_2 &= \text{Linear Stability Error} \label{eq: SBP_stability}
    \end{align}
\end{subequations}
These show that for small values of $\varepsilon$ the linear stability errors scale with $\varepsilon^2$ and conservation errors scale with $\varepsilon^6$. The loss of linear stability this shows is consistent with the rising value of maximum $\Im{(\hat{c}^\prime(\hat{k}))}$ shown in \cref{fig:diss_max}


To gain further insight, we pursue the combined analysis of \citet{Alhawwary2018}. Here we will use the diagonalisation of $\mathbf{Q}$:
\begin{equation}
    \mathbf{Q} = \mathbf{W\Lambda W}^{-1},
\end{equation}
which is valid in 1D provided the RBF centres and solution points are unique.

If we take the initial condition as being:
\begin{equation}
    \boldsymbol{u}_j(t=0) =\boldsymbol{u}_{j,0} = \exp{\left(ik\left(x_j + \frac{h}{2}(\pmb{\xi}+1)\right) \right)} = \pmb{\beta}\mathbf{W}
\end{equation}
then exponential integration can be used to find the approximate solution at time $t$ as:
\begin{equation}
    \tb{u}_j(t) = \mathbf{W}\exp{(-it\mathbf{\Lambda})}\mathbf{W}^{-1}\boldsymbol{u}_{j,0}
\end{equation}
Following \citet{Alhawwary2018}, we may then define the combined amplification factor as:
\begin{equation}
    G(k) = \frac{\|\tb{u}_{j}(t)\|_{L^2(\hat{\Omega})}}{\|\mathbf{u}_{j}(t)\|_{L^2(\hat{\Omega})}}
\end{equation}
and the phase difference as:
\begin{equation}
    \Delta\phi(k) = t - \angle\left(\int_{\hat{\Omega}} \tilde{u}_j(t)\overline{u}_j(t) \mathrm{d}\xi\right)
\end{equation}
where $\overline{z}$ and $\angle(z)$ are the complex conjugate and argument, respectively, of the complex number $z$.

\begin{figure}[tbhp]
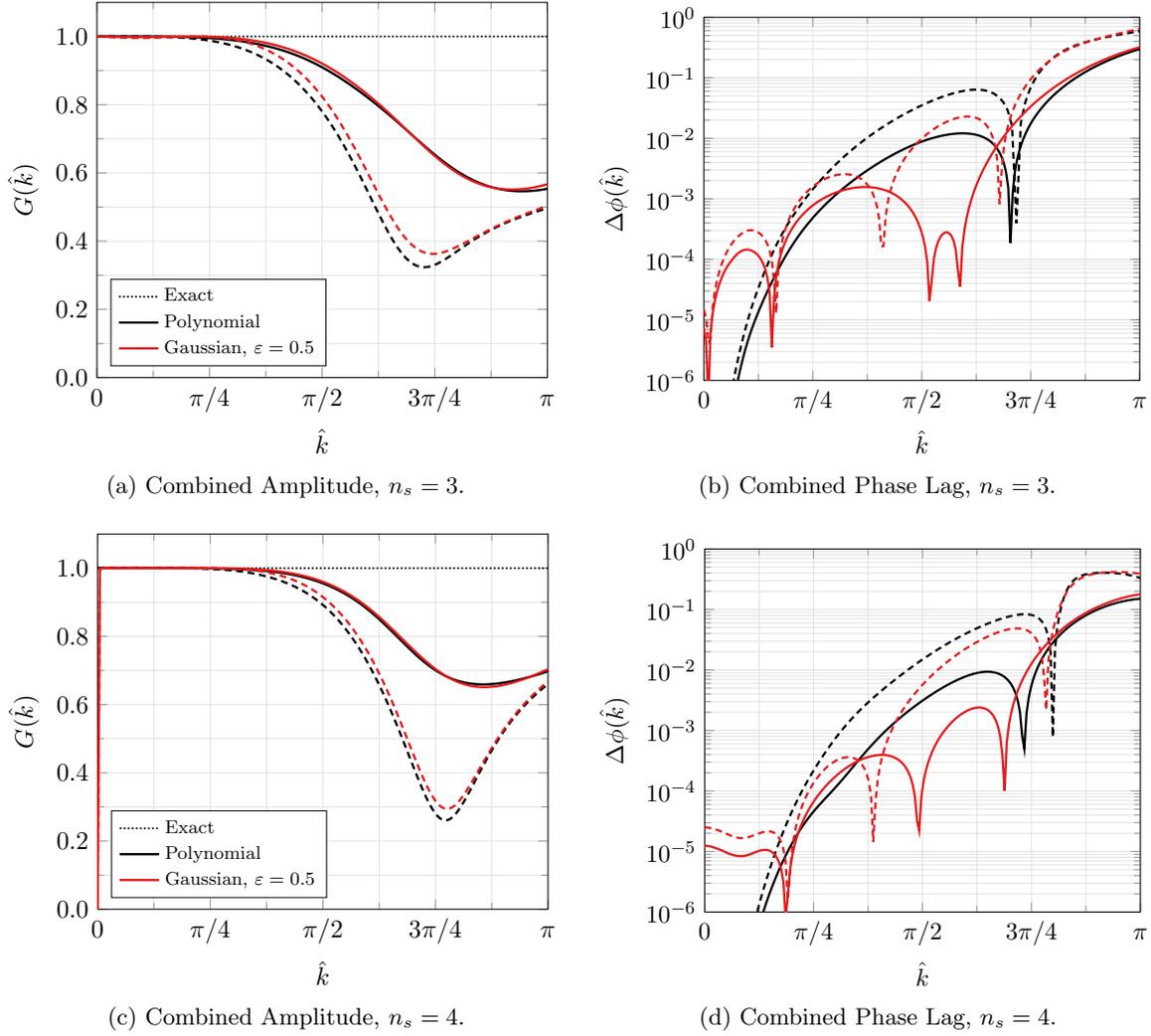

    \centering
    \subfloat[Combined Amplitude, $n_s=3$.]{\label{fig:cd_amp3}\adjustbox{width=0.49\linewidth, valign=b}{\input{Figs/tikz/n3_combined_fourier_amp}}}
    ~
    \subfloat[Combined Phase Lag, $n_s=3$.]{\label{fig:cd_phase3}\adjustbox{width=0.49\linewidth, valign=b}{\input{Figs/tikz/n3_combined_fourier_phase}}}
    \newline
    \subfloat[Combined Amplitude, $n_s=4$.]{\label{fig:cd_amp4}\adjustbox{width=0.49\linewidth, valign=b}{\input{Figs/tikz/n4_combined_fourier_amp}}}
    ~
    \subfloat[Combined Phase Lag, $n_s=4$.]{\label{fig:cd_phase4}\adjustbox{width=0.49\linewidth, valign=b}{\input{Figs/tikz/n4_combined_fourier_phase}}}
    \caption{\label{fig:cd}Combined semi-discrete Fourier analysis (solid $t=1$, dashed $t=2$).}
\end{figure}

\cref{fig:cd} presents the results of the amplification factor and phase lag at two instants in time, $t=1$ and $t=2$. These amplification factors show that the RBF bases produce lower dissipation, and that this difference becomes larger for longer integration times. This `holistic' approach confirms what was seen in the Fourier analysis results presented earlier. The holistic phase lag results present a more complex picture. Generally, the RBF basis is shown to have lower phase lag than the polynomial basis; however, at low wavenumbers for $n_s=4$, \cref{fig:cd_phase4}, a small constant phase error is observed, which could be due to the increasing numerical sensitivity of the RBFs at higher degrees. Numerical investigations will be required to determine what impact, if any, this has on the ultimate behaviour of the scheme.

\section{Numerical Experiments}\label{sec:numerical}

As well as the Fourier-inspired analysis of the equations, the effect of the RBF functional approximation was also investigated using numerical experiments of two kinds --- first, linear problems, typified by the linear advection--diffusion equation, and secondly nonlinear problems, using Burgers' equation.

\subsection{Linear Advection-Diffusion}\label{sec:linearadvdiff}
The basic one dimensional linear advection--diffusion equation is given by:
\begin{equation}\label{eq:linearAdvection}
    \px{u}{t} + a\px{u}{x} = \mu\pxi{2}{u}{x}, \quad \forall\; u(x,t):\Omega\times\mathbb{R}_+\mapsto\mathbb{R}
\end{equation}

Two sets of equation constants and initial conditions are used:
\begin{subequations}
    \begin{align}
        u(x,t=0) = \sin{x}&, \quad x\in\Omega=[-\pi,\pi), \quad \mathrm{with} \quad a=1, \mu=0\\
        u(x,t=0) = \exp{\left(-x^2\right)}&, \quad x\in\Omega=[-10,10), \quad \mathrm{with} \quad a=1, \mu=0.1
    \end{align}
\end{subequations}
which represent pure advection of a sine wave and advection--diffusion of a Gaussian bump. In both cases, the domain, $\Omega$, is taken as periodic.

To assess the resulting solutions, the point-mean $L^2$-norm of the error is used, with the approximate solution, $u^\delta(x,t)$, compared to the exact solution, $u(x,t)$. The following definition of the point-mean $L^2$ error norm is utilised:
\begin{equation}
    \|\mathcal{E}\|_2  = \sqrt{\frac{1}{N \times n_s} \sum^N_{j=1}\sum^{n_s}_{i=1}\left(\boldsymbol{u}^\delta_j(x_i) - u(x_{j,i})\right)^2}
\end{equation}
for a case with $N$ elements and with each element composed of $n_s$ solution points.

Along with the the $L^2$ norm, the $L^1$ and $L^\infty$ norms were also calculated. These were both found to behave similarly to the $L^2$ norm across all the cases tested, and so they are not further reported here.

\subsubsection{Layout of the Solution Points\label{ssec:solnpointlocs}}

\begin{figure}[tbhp]
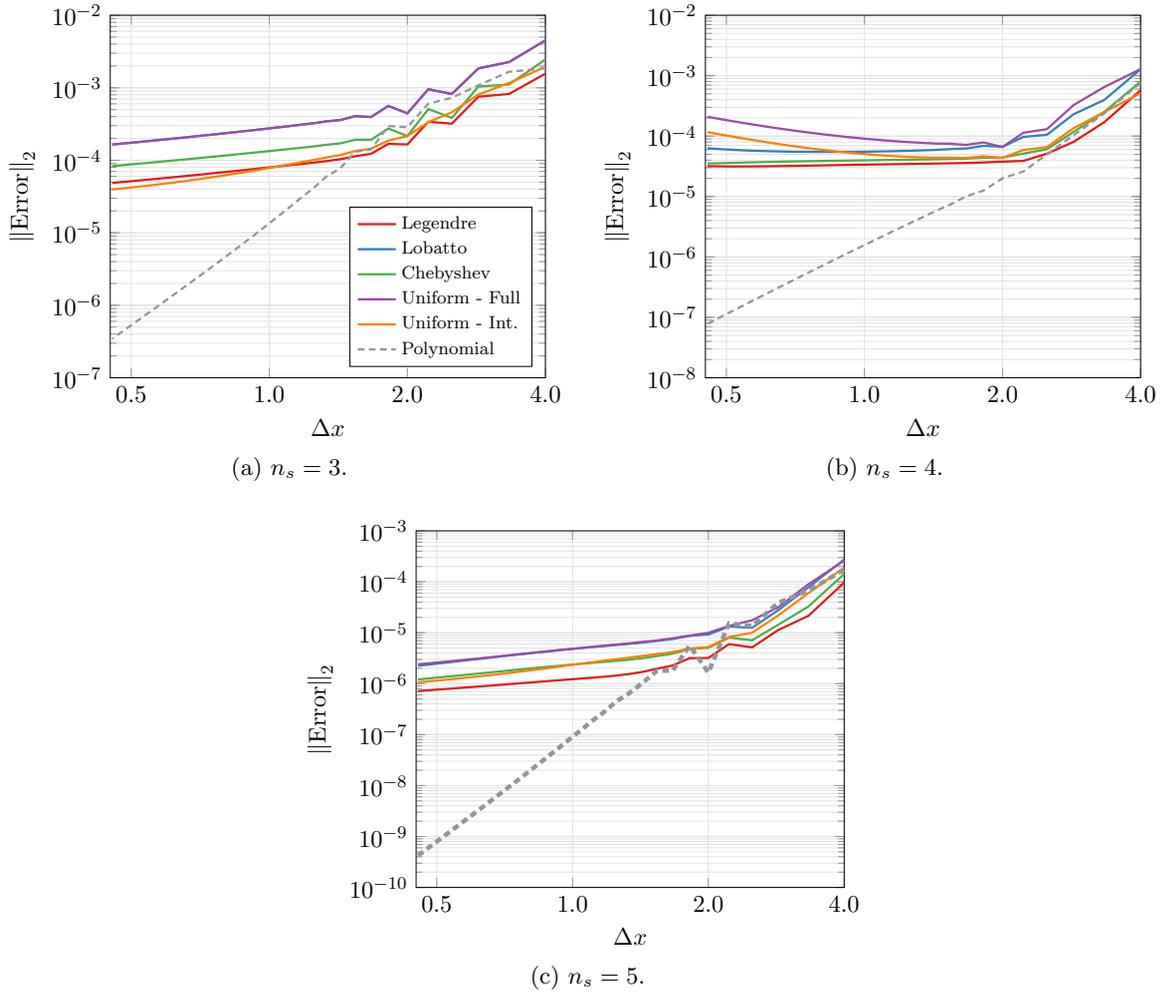

    \centering
    \subfloat[$n_s=3$.]{\label{fig:n3_gad_locations}\adjustbox{width=0.49\linewidth, valign=b}{\input{Figs/tikz/n3_gad_pointLocations}}}
    ~
    \subfloat[$n_s=4$.]{\label{fig:n4_gad_locations}\adjustbox{width=0.49\linewidth, valign=b}{\input{Figs/tikz/n4_gad_pointLocations}}}
    \newline
    \subfloat[$n_s=5$.]{\label{fig:n5_gad_locations}\adjustbox{width=0.49\linewidth, valign=b}{\input{Figs/tikz/n5_gad_pointLocations}}}
    \caption{\label{fig:gad_locations}Effect of point layout on $L^2$ error for linear advection--diffusion, $\varepsilon = 0.5$.}
\end{figure}

\begin{figure}[tbhp]
    \centering
    \subfloat[$n_s=3$.]{\label{fig:n3_sa_locations}\adjustbox{width=0.49\linewidth, valign=b}{\input{Figs/tikz/n3_sa_pointLocations}}}
    ~
    \subfloat[$n_s=4$.]{\label{fig:n4_sa_locations}\adjustbox{width=0.49\linewidth, valign=b}{\input{Figs/tikz/n4_sa_pointLocations}}}
    \newline
    \subfloat[$n_s=5$.]{\label{fig:n5_sa_locations}\adjustbox{width=0.49\linewidth, valign=b}{\input{Figs/tikz/n5_sa_pointLocations}}}
    \caption{\label{fig:sa_locations}Effect of point layout on $L^2$ error for linear advection, $\varepsilon = 0.5$.}
\end{figure}

\cref{fig:gad_locations,fig:sa_locations} show the performance of the RBF-FR method using different point locations for a range of mesh spacings. \cref{fig:gad_locations} shows behaviour of the advection--diffusion of the Gaussian bump, and~\cref{fig:sa_locations} shows the performance of pure advection of the sine wave, both at a modest shape parameter value of $\varepsilon=0.5$, and across three solution point numbers, of $n_s = 3$, $4$, and $5$. The functional centres were located at the solution points for all the cases. For comparison, a set of solutions computed with the polynomial bases using the Legendre points is also supplied. As expected, there is a strong dependency of the accuracy on the point layouts, with close to an order of magnitude of error separating the best performing points from the worst performing sets. Typically, the full uniform points and the Lobatto points perform least well, and, in general, the Legendre points perform best. For the advection--diffusion equation, there is some suggestion in \cref{fig:gad_locations} that the better performing solution point layouts are capable of outperforming the polynomial bases on the very coarsest meshes, but this does not hold true as the meshes are refined. The RBF-FR also seems to perform better, relative to the polynomial baseline, with $3$ and $5$ solution points, performing noticeably more poorly on finer meshes for the $n_s = 4$ case, for both the advection--diffusion and pure advection test cases. 

\begin{figure}[tbhp]
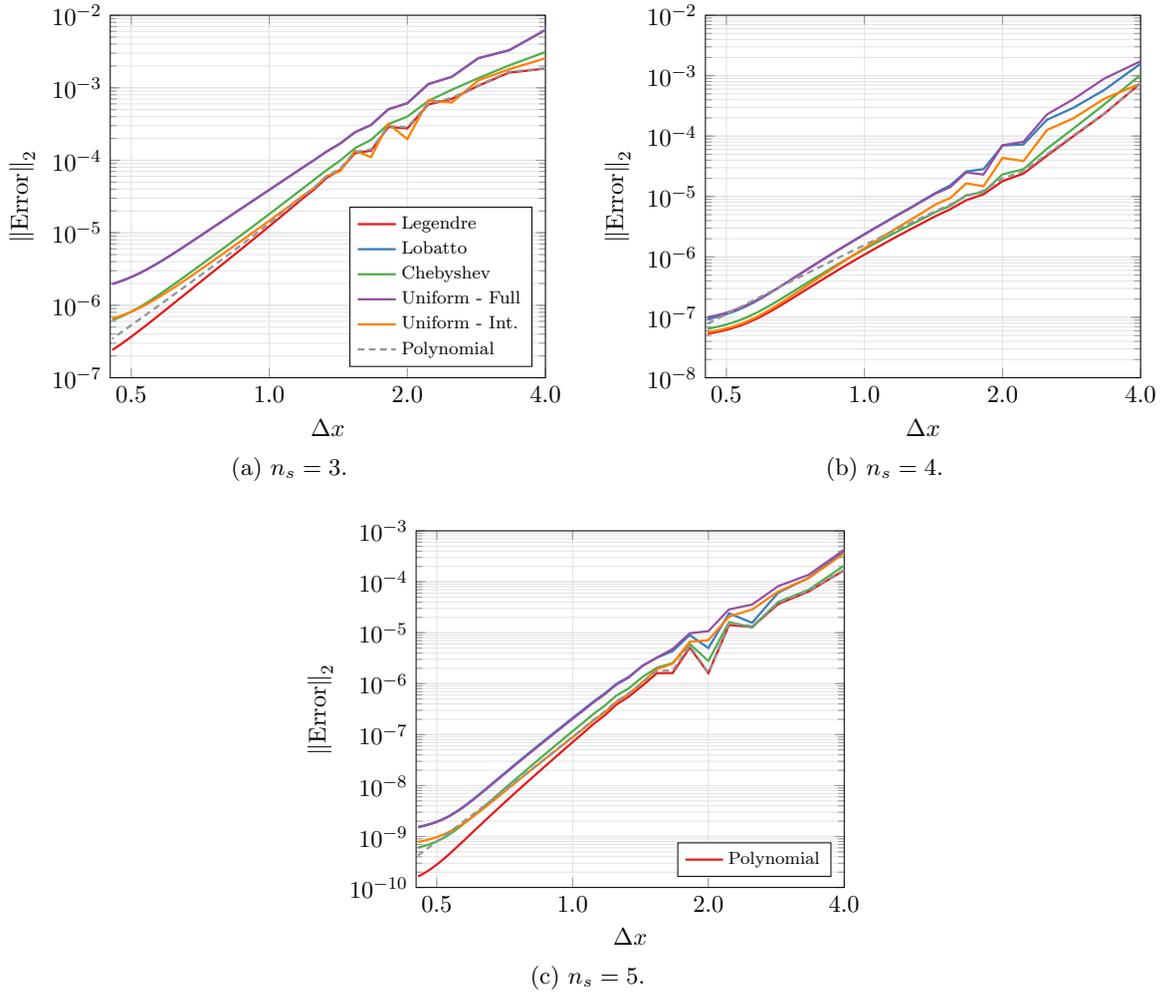

    \centering
    \subfloat[$n_s=3$.]{\label{fig:n3_gad_locations_01}\adjustbox{width=0.49\linewidth, valign=b}{\input{Figs/tikz/n3_gad_pointLocations_01}}}
    ~
    \subfloat[$n_s=4$.]{\label{fig:n4_gad_locations_01}\adjustbox{width=0.49\linewidth, valign=b}{\input{Figs/tikz/n4_gad_pointLocations_01}}}
    \newline
    \subfloat[$n_s=5$.]{\label{fig:n5_gad_locations_01}\adjustbox{width=0.49\linewidth, valign=b}{\input{Figs/tikz/n5_gad_pointLocations_01}}}
    \caption{\label{fig:gad_locations_01}Effect of point layout on $L^2$ error for linear advection--diffusion, $\varepsilon = 0.1$.}
\end{figure}

\cref{fig:gad_locations_01,fig:sa_locations_01} show the same set of runs as \cref{fig:gad_locations,fig:sa_locations}, but with the shape parameter of the radial bases reduced from $0.5$ to $0.1$. As would be expected, as the parameter is reduced towards zero, the performance of the fit follows the polynomial results more closely --- in the flat limit of $\varepsilon = 0$, the Gaussian RBF fit can be shown to be equivalent to a polynomial fit [REF]. In terms of the point locations, similar trends are observed with the wider Gaussian bases --- there is still close to an order of magnitude of difference between the best and worst performers. Again, the functions with colocated solution and flux points (the full uniform and Lobatto points) are seen to perform least effectively, and the Legendre points are again found to be the most effective of those tested. With the advection diffusion of the Gaussian bump, \cref{fig:gad_locations_01}, there is evidence that the best performing point locations outperform the polynomial FR scheme across the entire range tested. With the pure advection of the sine wave, \cref{fig:sa_locations_01}, the behaviour of RBF-FR is seen to be more complex. There are pronounced peaks and troughs in the error curves, which typically move the RBF-FR accuracy either side of the polynomial accuracy. Again, the general trend is that RBF-FR is able to outperform the polynomial method on the coarser meshes, but is overtaken as the meshes are refined.

\begin{figure}[tbhp]
    \centering
    \subfloat[$n_s=3$.]{\label{fig:n3_sa_locations_01}\adjustbox{width=0.49\linewidth, valign=b}{\input{Figs/tikz/n3_sa_pointLocations_01}}}
    ~
    \subfloat[$n_s=4$.]{\label{fig:n4_sa_locations_01}\adjustbox{width=0.49\linewidth, valign=b}{\input{Figs/tikz/n4_sa_pointLocations_01}}}
    \newline
    \subfloat[$n_s=5$.]{\label{fig:n5_sa_locations_01}\adjustbox{width=0.49\linewidth, valign=b}{\input{Figs/tikz/n5_sa_pointLocations_01}}}
    \caption{\label{fig:sa_locations_01}Effect of point layout on $L^2$ error for linear advection, $\varepsilon = 0.1$.}
\end{figure}

\subsubsection{Effect of the RBF Centre Locations}

In the previous sections, the functional centres of the RBFs have been placed at the same location as the solution points. Here, the effect of separating them is investigated.

\begin{figure}[tbhp]
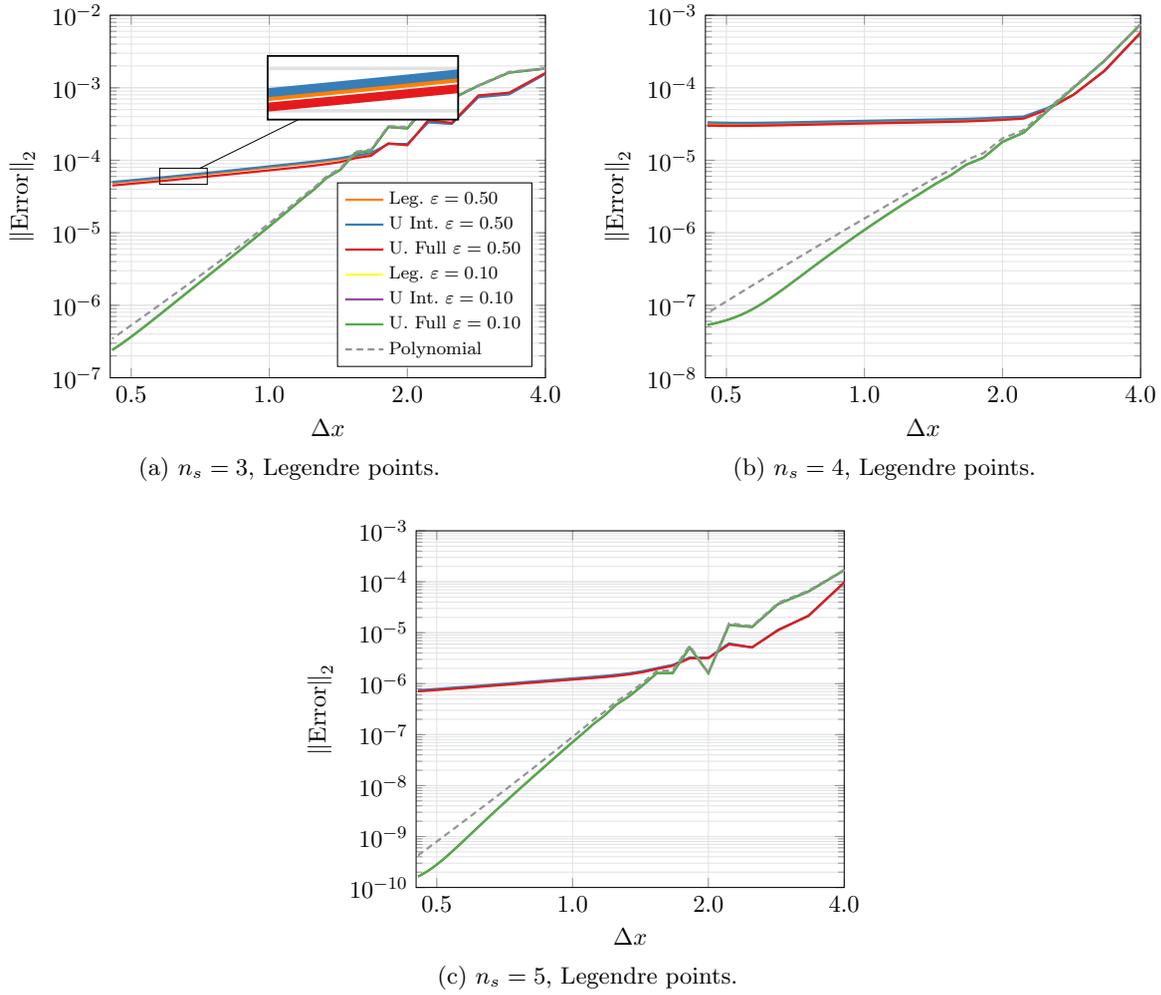

    \centering
    \subfloat[$n_s=3$, Legendre points.]{\label{fig:n3_gad_cenLocations}\adjustbox{width=0.49\linewidth, valign=b}{\input{Figs/tikz/n3_gad_cenLocations}}}
    ~
    \subfloat[$n_s=4$, Legendre points.]{\label{fig:n4_gad_cenLocations}\adjustbox{width=0.49\linewidth, valign=b}{\input{Figs/tikz/n4_gad_cenLocations}}}
    \newline
    \subfloat[$n_s=5$, Legendre points.]{\label{fig:n5_gad_cenLocations}\adjustbox{width=0.49\linewidth, valign=b}{\input{Figs/tikz/n5_gad_cenLocations}}}
    \caption{\label{fig:gad_cenLocations}Effect of varying the RBF centre locations on $L^2$ error for linear advection--diffusion.}
\end{figure}

\cref{fig:gad_cenLocations} shows these results, at two different values of shape parameter. The solution points are kept at the Legendre nodes, and the RBF centres are changed from the Legendre nodes to the full uniform and internal uniform distributions. From the plot, it appears that any influence made by moving the RBF centres is marginal, and, as might be expected, this influence reduces even further as the shape parameter decreases and the bases increasingly accurately approximate $\phi(r)=1$. In general, the results shown here for the advection--diffusion of the Gaussian bump suggest that there may be a small benefit in using the full uniform point set as the centres, but this effect is again marginal. It is possible that this small benefit is related to the effect of the functional centre locations on the condition number of the alternant matrix, and thus the numerical accuracy of calculating the operator matrices. This result, however, was not observed to be repeated for the pure advection of the sine wave, and, due to the marginal nature of any effect seen, further results are omitted here for brevity.

\subsubsection{Effect of the Shape Parameter\label{ssec:shapefunction}}

Given its importance to RBF interpolation, the effect of reducing this scaling parameter is explored through the two test cases. In all cases, the best performing point locations from the previous section were used --- the Legendre points --- and the functional centres were located at the solution points. 

\begin{figure}[tbhp]
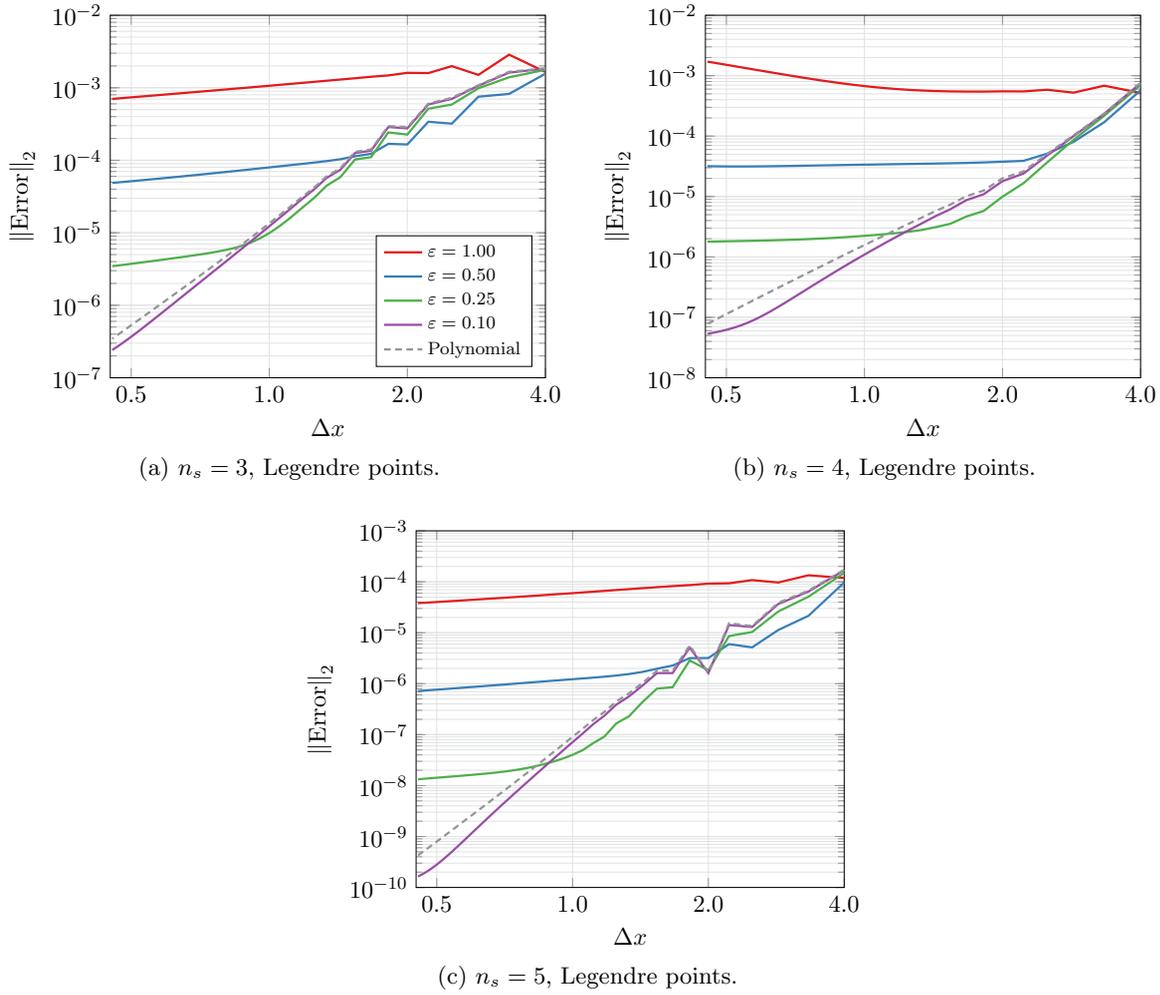

    \centering
    \subfloat[$n_s=3$, Legendre points.]{\label{fig:n3_gad_varEpsilon}\adjustbox{width=0.49\linewidth, valign=b}{\input{Figs/tikz/n3_gad_varEpsilon}}}
    ~
    \subfloat[$n_s=4$, Legendre points.]{\label{fig:n4_gad_varEpsilon}\adjustbox{width=0.49\linewidth, valign=b}{\input{Figs/tikz/n4_gad_varEpsilon}}}
    \newline
    \subfloat[$n_s=5$, Legendre points.]{\label{fig:n5_gad_varEpsilon}\adjustbox{width=0.49\linewidth, valign=b}{\input{Figs/tikz/n5_gad_varEpsilon}}}
    \caption{\label{fig:gad_varEpsilon}Effect of shape parameter, $\varepsilon$, on $L^2$ error for linear advection--diffusion.}
\end{figure}

\begin{figure}[tbhp]
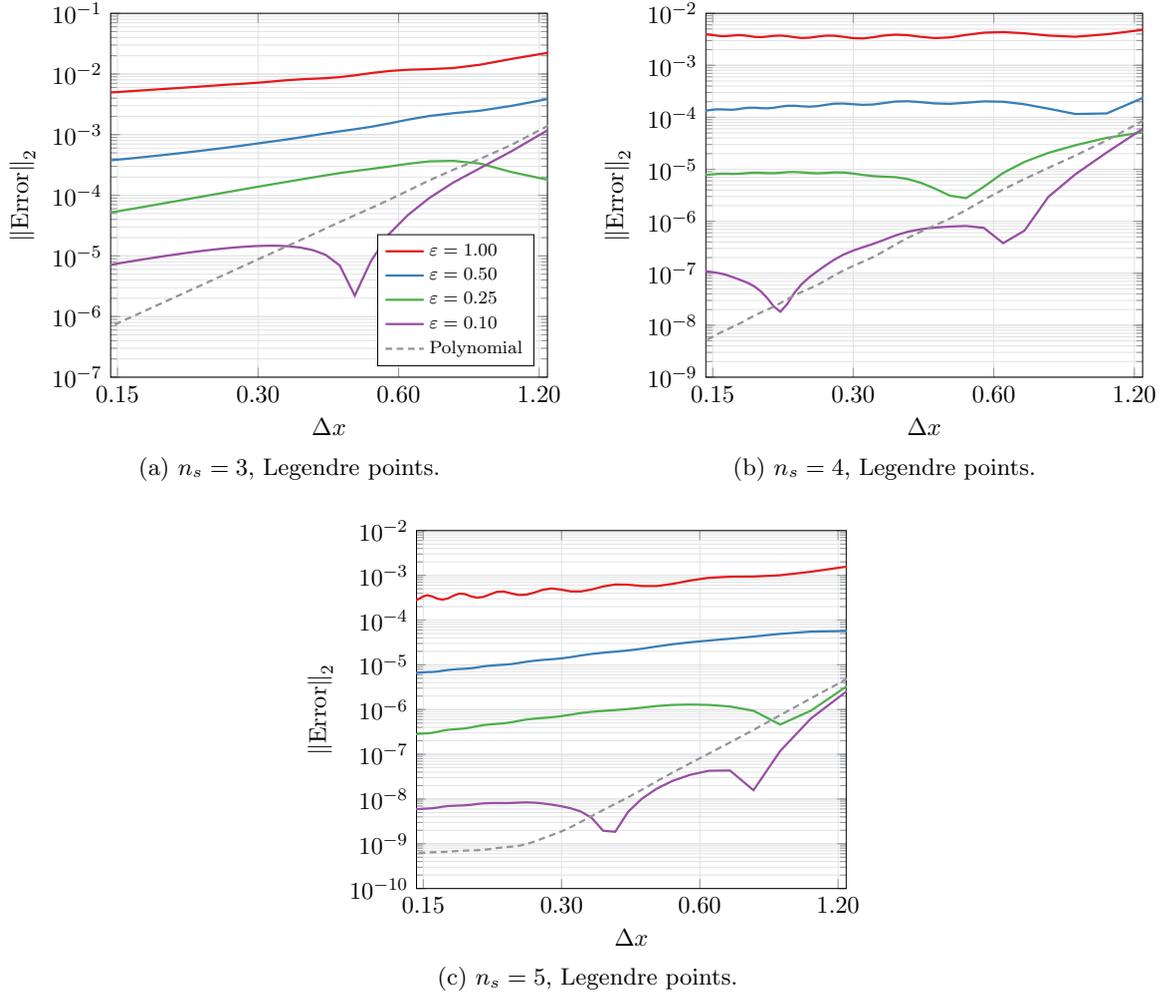

    \centering
    \subfloat[$n_s=3$, Legendre points.]{\label{fig:n3_sa_varEpsilon}\adjustbox{width=0.49\linewidth, valign=b}{\input{Figs/tikz/n3_sa_varEpsilon}}}
    ~
    \subfloat[$n_s=4$, Legendre points.]{\label{fig:n4_sa_varEpsilon}\adjustbox{width=0.49\linewidth, valign=b}{\input{Figs/tikz/n4_sa_varEpsilon}}}
    \newline
    \subfloat[$n_s=5$, Legendre points.]{\label{fig:n5_sa_varEpsilon}\adjustbox{width=0.49\linewidth, valign=b}{\input{Figs/tikz/n5_sa_varEpsilon}}}
    \caption{\label{fig:sa_varEpsilon}Effect of shape parameter, $\varepsilon$, on $L^2$ error for linear advection.}
\end{figure}

\cref{fig:gad_varEpsilon,fig:sa_varEpsilon} show the effect of reducing this parameter on the performance of RBF-FR, again with a polynomial FR curve for comparison, for $n_s = 3$, $4$, and $5$. Generally, the accuracy of the RBF-FR solution was found to improve as the shape parameter is reduced. \cref{fig:gad_varEpsilon}, for the bump advection--diffusion case, suggests that each of the shape parameter values are capable of outperforming the polynomial fit across some range of grid densities. Larger values of $\varepsilon$ tend to outperform the polynomial on extremely coarse grids, but as the grid is refined, the polynomial approach overtakes them. As the value of the shape parameter is reduced, the range of mesh densities at which RBF-FR outperforms polynomial FR moves towards finer meshes. By the time a value of $\varepsilon = 0.1$ is reached, the RBF-FR outperforms the polynomial across the entire range of mesh densities tested, with close to a $50\%$ reduction in the error on the finest mesh tested for $n_s = 5$. The other observable trend is that as the shape parameter is reduced, the performance becomes increasingly like that of the polynomial FR. 

\cref{fig:sa_varEpsilon} shows the results of the same range of tests, but for the pure advection of a sine wave. The RBF-FR performs much less competitively here, but similar trends can be observed --- higher values of shape parameter, so narrower bases, tend to outperform the polynomial scheme on more coarse meshes, and wider bases tend to outperform the polynomial FR on finer meshes. However, in all of these cases the polynomial FR outperforms RBF-FR for the finest meshes. Despite this, there are still large ranges of mesh densities in which the RBF-FR is highly competitive with the polynomial FR in terms of $L^2$ error norm, outperforming it by an order of magnitude in some cases. As seen with the advection diffusion of the Gaussian bump, as the shape parameter is reduced, the curves gradually trend towards the baseline polynomial error curve. With the baseline Gaussian RBF, it was not possible to use shape parameters significantly smaller than $\varepsilon = 0.1$, as the condition number of the system of equations requiring solution became too large, and numerical errors swamped the solution.

\subsubsection{Small Values of Shape Parameter\label{subsec:smallShapeParameters}}

Using RBF-Direct to compute the necessary operators limits the minimum possible value of $\varepsilon$. This issue can be sidestepped by making use of the RBF-GA method to compute the operators.

\begin{figure}[tbhp]
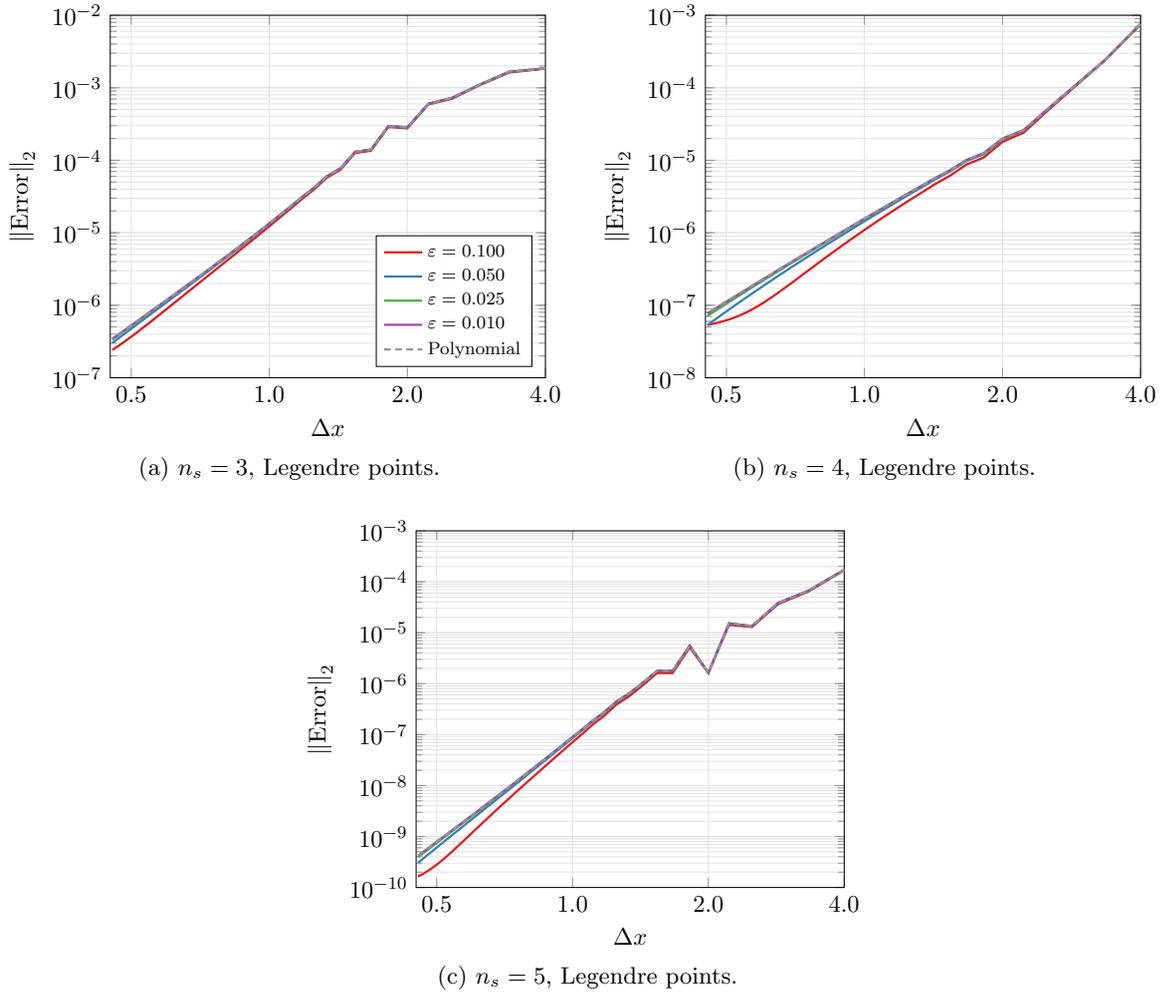

    \centering
    \subfloat[$n_s=3$, Legendre points.]{\label{fig:n3_gad_smallEpsilonc}\adjustbox{width=0.49\linewidth, valign=b}{\input{Figs/tikz/n3_gad_smallEpsilon}}}
    ~
    \subfloat[$n_s=4$, Legendre points.]{\label{fig:n4_gad_smallEpsilon}\adjustbox{width=0.49\linewidth, valign=b}{\input{Figs/tikz/n4_gad_smallEpsilon}}}
    \newline
    \subfloat[$n_s=5$, Legendre points.]{\label{fig:n5_gad_smallEpsilon}\adjustbox{width=0.49\linewidth, valign=b}{\input{Figs/tikz/n5_gad_smallEpsilon}}}
    \caption{\label{fig:gad_smallEpsilon}Effect of very small shape parameters on $L^2$ error for linear advection--diffusion.}
\end{figure}

\cref{fig:gad_smallEpsilon,fig:sa_smallEpsilon} show the effect of using smaller values of $\varepsilon$ has on the performance of the scheme, made possible by RBF-GA. To test for consistency between the two methods, a value of $\varepsilon=0.1$ is repeated from the RBF-Direct approach, and is found to give the same results as the RBF-GA method for all cases --- this is as expected, given that both have been deliberately constructed to span the same functional approximation basis space.

At lower values of shape parameter, the sweep continues to trend towards the polynomial result, and pushes the range of mesh densities at which RBF-FR outperforms the polynomial scheme further to the left. For the pure advection of the sine wave, where RBF-FR has previously performed poorly relative to the polynomial scheme, very low values of shape parameter tend to produce more competitive results, with values of $\varepsilon \leq 0.025$ outperforming polynomial approximation across the entire range tested. For $n_s=5$, all values of $\varepsilon \leq 0.05$ outperform or match the polynomial scheme.

The use of RBF-GA to calculate the bases gives RBF-FR the ability to outperform polynomial FR on both test cases at no additional runtime cost, while retaining the geometric flexibility of both RBFs and FR in higher dimensions. 

\begin{figure}[tbhp]
    \centering
    \subfloat[$n_s=3$, Legendre points.]{\label{fig:n3_sa_smallEpsilon}\adjustbox{width=0.49\linewidth, valign=b}{\input{Figs/tikz/n3_sa_smallEpsilon}}}
    ~
    \subfloat[$n_s=4$, Legendre points.]{\label{fig:n4_sa_smallEpsilon}\adjustbox{width=0.49\linewidth, valign=b}{\input{Figs/tikz/n4_sa_smallEpsilon}}}
    \newline
    \subfloat[$n_s=5$, Legendre points.]{\label{fig:n5_sa_smallEpsilon}\adjustbox{width=0.49\linewidth, valign=b}{\input{Figs/tikz/n5_sa_smallEpsilon}}}
    \caption{\label{fig:sa_smallEpsilon}Effect of very small shape parameters on $L^2$ error for linear advection.}
\end{figure}

\subsection{\label{ssec:bt}Viscous Burgers' Equation}
    The viscous Burgers' equation represents a meaningful step towards the full Navier--Stokes equations. Burgers' equation is nonlinear, and is able to display turbulence-like phenomena, making it an insightful numerical test case as it exhibits integral, inertial, and dissipation ranges. Several other works have studied the phenomenon of Burgers' turbulence, notably those of \citet{San2016,Alhawwary2018,Frisch2013}. The primary relevant difference compared to the Navier--Stokes equations is that the energy cascade of Burgers' turbulence follows a $k^{-2}$ profile, compared to the famous $k^{-5/3}$ relation seen in Navier--Stokes.
    
    For the viscous Burgers' equation:
    \begin{equation}\label{eq:burgers}
        \px{u}{t} + \half\px{u^2}{x} = \mu\pxi{2}{u}{x} \quad \text{for} \quad u(x,t) : \Omega\times\mathbb{R}_+ \mapsto \mathbb{R}
    \end{equation}
    the initial condition proposed by \citet{San2016} is used, where energy is distributed among modes according to the spectrum given by:
	\begin{equation}
	    E_0(k) = \frac{2}{3\sqrt{\pi}\rho}(k\rho)^4\exp{\left(-(k\rho)^2\right)}, \quad \mathrm{where} \quad \rho = \frac{1}{10}
    \end{equation}
    Use of these values gives the most energised mode as corresponding to $k=13$. For a periodic domain of $\Omega=[0,2\pi)$, the solution field may then be reconstructed from the spectrum as:
    \begin{equation}
	    u_0(x) = \sum^{k_\mathrm{max}}_{k=0}\sqrt{2E_0(k)}\cos{(kx + 2\pi\Psi(k))}
    \end{equation}
    where $k_\mathrm{max} = 2048$ is found to be sufficient for double-precision floating-point accuracy. The term $\Psi(k)\in[0,1)$ is used to produce a random phase angle for each wavenumber, and for reproducibility the Mersenne twister algorithm~\citep{Matsumoto1998} is used.
    
    The quantity of interest is the frequency response of the system, and to understand this the resulting energy spectrum at $t=0.1$ is studied. Time integration is performed using an RK4 scheme and, to enable discrete Fourier transforms, data is interpolated via the underlying basis to a set of equispaced points. For a uniform one-dimensional mesh with spacing $h$, uniform reference points are located at $\hat{x}\in\{-1+\hat{h}, -1+3\hat{h},\dots, -1+(2n-1)\hat{h}\}$, for $\hat{h}=2/(2n)$. The discrete Fourier transform results were averaged over 1000 runs, with the same 1000 random seeds used for the different initial states tested. Using the same 1000 seeds is important for objective comparison, as without this the behaviour at high frequencies can appear to be quite different between tests.

    \begin{figure}[tbhp]
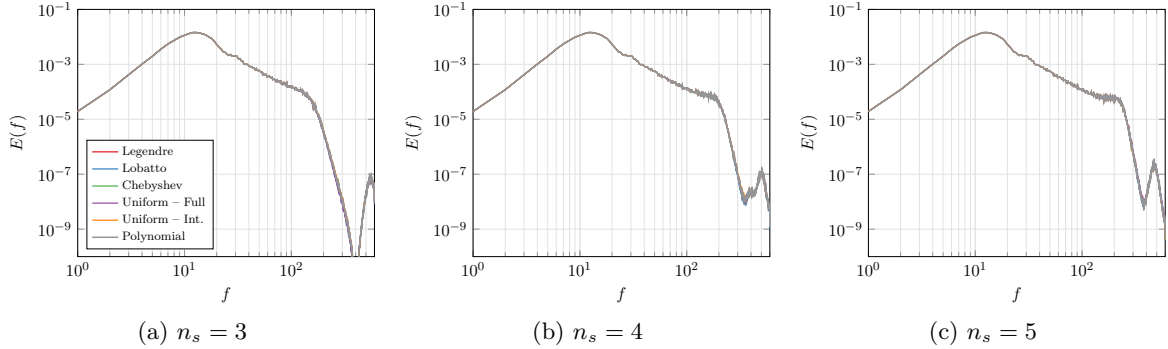
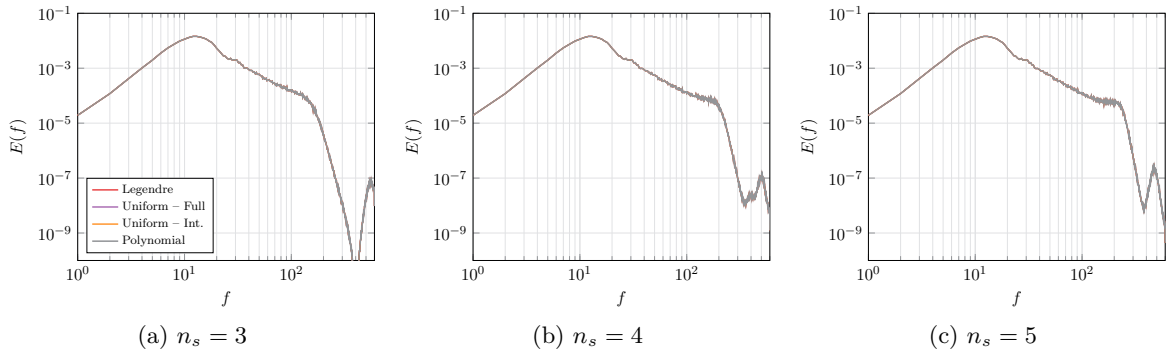

        \centering
        \subfloat[$n_s=3$]{\label{fig:bt_pnt_03}\adjustbox{width=0.32\linewidth, valign=b}{\input{Figs/tikz/bt_points_n03}}}
        ~
        \subfloat[$n_s=4$]{\label{fig:bt_pnt_04}\adjustbox{width=0.32\linewidth, valign=b}{\input{Figs/tikz/bt_points_n04}}}
        ~
        \subfloat[$n_s=5$]{\label{fig:bt_pnt_05}\adjustbox{width=0.32\linewidth, valign=b}{\input{Figs/tikz/bt_points_n05}}}
        \caption{\label{fig:bt_pnt}Effect of changing solution and collocated centre locations for Gaussian RBF with $\varepsilon=0.01$.}
    \end{figure}
    
    Investigating first the effect of varying the solution points (with colocated RBF centres), it can be seen from \cref{fig:bt_pnt} that there seems to be no significant impact on the energy spectra of the Burgers' turbulence for the point layouts of \cref{fig:pointLayouts}. In particular, it is known that the turbulent bottleneck~\citep{Fernandez2019} is caused by higher order numerical diffusion terms resulting from the discretisation. Therefore, if the scale of these terms was meaningfully changed by using either an RBF approximation space or by the point locations, it would show as a change in the bottleneck effect. The bottleneck is clearly visible in the approximate range $100<f<250$, and is unchanged here.
    
    These tests were repeated, but fixing the solution points at the Legendre points and varying the locations of the RBF centres. The results of these tests are presented in \cref{fig:bt_ctr}. Once again, no significant variation in the Burgers' turbulence spectra is seen. 
    
    \begin{figure}[tbhp]
        \centering
        \subfloat[$n_s=3$]{\label{fig:bt_ctr_03}\adjustbox{width=0.32\linewidth, valign=b}{\input{Figs/tikz/bt_centre_n03}}}
        ~
        \subfloat[$n_s=4$]{\label{fig:bt_ctr_04}\adjustbox{width=0.32\linewidth, valign=b}{\input{Figs/tikz/bt_centre_n04}}}
        ~
        \subfloat[$n_s=5$]{\label{fig:bt_ctr_05}\adjustbox{width=0.32\linewidth, valign=b}{\input{Figs/tikz/bt_centre_n05}}}
        \caption{\label{fig:bt_ctr}Effect of changing solution and collocated centre locations for Gaussian RBF with $\varepsilon=0.01$.}
    \end{figure}
    
    \begin{figure}[tbhp]
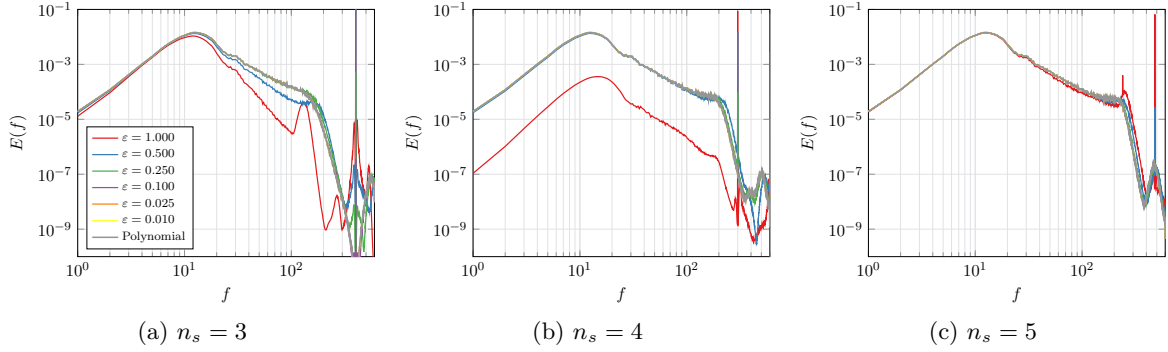

        \centering
        \subfloat[$n_s=3$]{\label{fig:bt_eps_03}\adjustbox{width=0.32\linewidth, valign=b}{\input{Figs/tikz/bt_eps_n03}}}
        ~
        \subfloat[$n_s=4$]{\label{fig:bt_eps_04}\adjustbox{width=0.32\linewidth, valign=b}{\input{Figs/tikz/bt_eps_n04}}}
        ~
        \subfloat[$n_s=5$]{\label{fig:bt_eps_05}\adjustbox{width=0.32\linewidth, valign=b}{\input{Figs/tikz/bt_eps_n05}}}
        \caption{\label{fig:bt_eps}Effect of varying $\varepsilon$ for Legendre solution points and internal uniform centres.}
    \end{figure}
    
    The final test undertaken with Burgers' turbulence was to again vary the RBF shape parameter, $\varepsilon$. The results shown in \cref{fig:bt_eps} were produced with the Legendre solution points and the internal uniform functional centres. This test clearly shows that $\varepsilon$ can have a profound impact on the behaviour of the method. As $\varepsilon$ is increased, these results appear to show that the method becomes unstable at high frequencies, with the width of the peak caused by this instability spreading. This is in agreement with the dissipation results of the Fourier analysis. Furthermore, as $\varepsilon$ is increased the location of the bottleneck is shifted towards higher frequencies. This is consistent with an overall reduction in the dissipation, and is again in agreement with the results of the Fourier analysis presented earlier. 
\section{Conclusions}\label{sec:conclusions}

In this paper, the Gaussian RBF has been used to replace polynomials as the method of approximating functions within elements as part of the broader FR framework for solving PDEs The effect that this has on the performance of the simulations is investigated both analytically and numerically.

The shape parameter, $\varepsilon$, is found to be of vital importance. Lower values of shape parameter trend towards the polynomial performance, but there tends to be a range of mesh densities at which the RBF-FR approach is capable of outperforming the polynomial, for linear problems. This range of mesh densities becomes finer as the shape parameter is reduced.

Although directly fitting RBFs to scattered points becomes untenable for very low values of $\varepsilon$, a method was discussed to avoid this. Although this makes the operators more expensive to compute, the nature of flux reconstruction means that these operator matrices are calculated \emph{a priori}, and do not need to be recomputed for each run, meaning there is no increase in simulation run time. Very low values of shape parameter are proven to make RBF-FR effective for both linear and nonlinear problems.

The position of the solution points is also demonstrated to have a considerable impact on the performance of the simulation, with the Legendre points being found to be the most effective for linear problems. This difference was less pronounced for the nonlinear numerical experiments. Although this paper is restricted to studying the effects in one dimension, the increased geometrical flexibility of RBFs make them potentially attractive in higher dimensions. Future work will explore the potential of this geometric flexibility in improving the performance of FR.

The choice of functional centres for the RBFs was found to have a very small influence on the performance --- particularly as the shape parameter was reduced and the RBFs became flat.

\section*{Acknowledgements}\label{sec:ack}

\section*{Data Availability}
The data that informed the findings of this study are available from the corresponding author upon reasonable request.

\bibliographystyle{plainnat}
\bibliography{reference}


\begin{appendices}
\end{appendices}


\end{document}